\documentclass[11pt,reqno]{amsart}
\usepackage[margin=1.05in]{geometry}
\usepackage{amsmath,amssymb,amsthm,mathtools,booktabs,array}
\usepackage{microtype}
\usepackage[colorlinks=true,linkcolor=blue,citecolor=blue,urlcolor=blue]{hyperref}

\usepackage{tikz}
\usetikzlibrary{arrows.meta,decorations.pathreplacing,calc}
\usepackage[T1]{fontenc}
\usepackage{mathrsfs}

\usepackage[T1]{fontenc}
\usepackage{lmodern}

\newcommand{\E}{\mathbb E}
\newcommand{\F}{\mathbb F}

\newcommand{\Z}{\mathbb Z}
\newcommand{\C}{\mathbb C}

\newcommand{\Lam}{\Lambda}
\newcommand{\lc}{\operatorname{lc}}
\newcommand{\spn}{\operatorname{span}}
\newcommand{\cE}{\mathcal E}
\newcommand{\cC}{\mathcal C}
\newcommand{\cB}{\mathcal B}

\newcommand{\wh}{\widehat}
\newtheorem{theorem}{Theorem}[section]
\newtheorem{proposition}[theorem]{Proposition}
\newtheorem{lemma}[theorem]{Lemma}
\newtheorem{corollary}[theorem]{Corollary}

\theoremstyle{definition}
\newtheorem{definition}[theorem]{Definition}
\theoremstyle{remark}
\newtheorem{remark}[theorem]{Remark}

\numberwithin{equation}{section}
\title{On the multidimensional polynomial Szemer\'{e}di theorem}
\author[D. He]{Danqing He}
\address{Danqing He, School of Mathematical Sciences,
Fudan University, People's Republic of China}
\email{hedanqing@fudan.edu.cn}
\thanks{D. He was supported by  National Key R$\&$D Program of China (No. 2021YFA1002500), NNSF of China (No. 12322105), and  the New Cornerstone Science Foundation.}
\author[X. Shao]{Xuancheng Shao}
\address{Department of Mathematics, University of Kentucky\\
Lexington, KY 40506\\
USA}
\email{xuancheng.shao@uky.edu}
\thanks{XS was supported by NSF grant DMS-2452462.}
\begin{document}
\begin{abstract}

We study the multidimensional polynomial Szemer\'{e}di theorem over
prime fields. Let $P_1,\dots,P_m\in\F_p[t]^k$ be block independent
vector polynomials of degree at most $d$, with $P_i(0)=0$ for every
$i$. We prove that there is a constant $c=c(k,m,d)>0$ such that,
for every sufficiently large prime $p$, every set $A\subset\F_p^k$
with $|A|\ge p^{k-c}$ contains a configuration
$x,x+P_1(t),\dots,x+P_m(t)$ with $t\in\F_p\setminus\{0\}$.

\end{abstract}
\maketitle
\tableofcontents

\section{Introduction}\label{sec:intro}
Let $P_1,\ldots,P_m\in\F_p[t]^k$ be vector polynomials having zero constant terms over the finite field $\F_p$. We consider the $k$-dimensional polynomial  progression
\[
 x,\ \ x+P_1(t),\ \ \cdots,\ \ x+P_m(t),
\]
where $x\in\F_p^k$ and $t\in\F_p$.
The purpose of this paper is to obtain an asymptotic formula for the number of such polynomial progressions under a natural independence condition on $P_1,\cdots,P_m$. 

\subsection{Main results}
To state our main results, we need two different spaces associated with $P_1,\cdots,P_m$. For $1 \leq i \leq m$, write $P_i = (P_{i,1},\cdots,P_{i,k})$ where each $P_{i,j} \in \F_p[t]$ is a polynomial, and write $P_i(t) = \sum_{n\geq 1}a_{i,n}t^n$ where the coefficients $a_{i,n} \in \F_p^k$. Define
\begin{equation}\label{eq:VW}
V_i = V(P_i) = \spn\{a_{i,n}:n\ge1\} \subset \F_p^k, \ \ 
 W_i = W(P_i) =\spn\{P_{i,1},\ldots,P_{i,k}\} \subset \F_p[t],
\end{equation}
where the spans are over $\F_p$.
We always have $\dim V_i = \dim W_i$.  

For a subspace $V\subset \F_p^k$ and a function $f: \F_p^k\rightarrow \C$, write $\E(f\vert V): \F_p^k\rightarrow \C$ for the average over cosets of $V$:
\[
 \E(f\vert V)(x)=\E_{v\in V}f(x+v).
\]
Thus $\E(f\vert V)$ is constant on each coset of $V$. For $P = (P_1,\cdots,P_m)$ and functions $f_0,\ldots,f_m: \F_p^k\rightarrow\C$, define the counting operator
$$
 \Lam_P(f_0,\ldots,f_m)
 =\E_{x\in \F_p^k} \E_{t\in\F_p}f_0(x)\prod_{i=1}^m f_i(x+P_i(t)).
$$ 

\begin{theorem}\label{thm:main}
Let $P = (P_1,\cdots,P_m)$, where $P_1,\ldots,P_m\in\F_p[t]^k$ with degree at most $d$ and $P_i(0)=0$ for each $1 \leq i \leq m$. Let $V_i,W_i$ be defined as in \eqref{eq:VW}. Suppose that
\begin{equation}\label{eq:independence}
 W_1+\cdots+W_m=W_1\oplus\cdots\oplus W_m.
\end{equation}
Then for all $1$-bounded functions $f_0,\ldots,f_m:\F_p^k\to\C$, we have
$$
 \Lam_P(f_0,\ldots,f_m) = \E_{x\in\F_p^k}f_0(x)\prod_{i=1}^m\E(f_i\vert V_i)(x) + O_{k,m,d}( p^{-c})
$$
for some constant $c = c(k,m,d) > 0$.
\end{theorem}

The independence assumption is necessary for a uniform formula of this form. If~\eqref{eq:independence} fails, there are rational linear forms $\ell_i$ such that
\[
 \sum_{i=1}^m\ell_i(P_i(t))=0,
 \qquad \ell_j(P_j(t))\ne0\quad\text{for at least one }j.
\]
Reduce these forms modulo a sufficiently large prime and take
\[
 f_i(x)=\exp\bigl(2\pi i\ell_i(x)/p\bigr)\quad(1\le i\le m),
 \qquad f_0(x)=\overline{\prod_{i=1}^m f_i(x)}.
\]
The polynomial relation makes the count equal to one. On the other hand, $\ell_j$ is nonzero on $V_j$, so $\E(f_j\vert V_j)=0$ and the proposed main term vanishes. Dependencies within a single block create no such obstruction; they are accounted for by its coefficient space.

Theorem \ref{thm:main} generalizes a result of Kuca \cite{Kuca} who  studied the case when $P_i(t) = q_i(t)v_i$ for some scalar polynomial $q_i(t)$ and nonzero vector $v_i \in \F_p^k$, under the assumption that $q_1,\cdots,q_m$ are linearly independent. In this case, each $V_i$ is a line spanned by $v_i$.

Under the stricter assumption that all the $km$ scalar polynomials $P_{i,j}$ in the components of $P_1,\cdots,P_m$ are linearly independent, one can show that $\dim V_i = \dim W_i = k$ and hence $V_i = \F_p^k$ for each $i$. This leads to the following corollary, which can be viewed as a multi-dimensional version of Peluse's result \cite{Peluse}.

\begin{corollary}\label{cor:full-rank}
Let $P = (P_1,\cdots,P_m)$, where $P_i = (P_{i,1},\cdots,P_{i,k})\in\F_p[t]^k$  with degree at most $d$ and $P_i(0)=0$ for each $1 \leq i \leq m$.
Suppose that the $km$ scalar polynomials $P_{i,j}$ are linearly independent over $\F_p$. Then for all $1$-bounded functions $f_0,\ldots,f_m:\F_p^k\to\C$, we have
\[
 \Lam_P(f_0,\ldots,f_m)=\prod_{i=0}^m\E f_i+O_{k,m,d}(p^{-c})
\]
for some constant $c = c(k,m,d)> 0$.
\end{corollary}

For example, consider
\[
 P_1(t)=(t^4,t^4+t),\qquad
 P_2(t)=(t^4+t^2,t^4+t^3).
\]
The four scalar coordinates are linearly independent, as is seen by comparing the coefficients of $t,t^2,t^3,t^4$. Thus Corollary~\ref{cor:full-rank} applies, although all four coordinates have the same degree.

By applying Theorem \ref{thm:main} to $f_0 = \cdots = f_m = 1_A$ for a subset $A \subset \F_p^k$, we can deduce the following quantitative multi-dimensional polynomial Szemer\'edi theorem under the assumption \eqref{eq:independence}.

\begin{corollary}\label{cor:szemeredi}
Let $P = (P_1,\cdots,P_m)$, where $P_1,\ldots,P_m\in\F_p[t]^k$ with degree at most $d$ and $P_i(0)=0$ for each $1 \leq i \leq m$. Let $W_i$ be defined as in \eqref{eq:VW}, and suppose that \eqref{eq:independence} holds. Then for any subset $A \subset \F_p^k$ with $|A| = \delta p^k$, we have
$$
\#\{(x,t) \in \F_p^k\times \F_p: x, x+P_1(t),\cdots,x+P_m(t) \in A\} \geq (\delta^{m+1} + O_{k,m,d}(p^{-c}))p^{k+1}
$$
for some constant $c = c(k,m,d) > 0$.
\end{corollary}

Here we use that $\delta^{m+1}\le \E_{x\in\F_p^k}f_0(x)\prod_{i=1}^m\E(f_i\vert V_i)(x)$, which follows  from $\E(f/\E(f|V_i))\le 1$ and H\"older's inequality.

\subsection{Relation to prior works}

The polynomial Szemer\'edi theorem of Bergelson and Leibman~\cite{BL} guarantees general polynomial configurations with zero constant terms in positive-density subsets of $\Z^k$, without an independence assumption on the polynomials. Bergelson, Leibman and McCutcheon~\cite{BLM} extended polynomial recurrence to modules over integral domains. These results provide the qualitative background for our problem. Under the independence assumption~\eqref{eq:independence}, we seek the more precise conclusion of a counting formula with a power-saving error.

Bergelson, Host, McCutcheon and Parreau \cite{BHMP} established a qualitative Szemer\'edi Theorem for vector-valued polynomials of one parameter in $\mathbb Z^k$, and Corollary~\ref{cor:szemeredi} can be regarded as a quantitative version of their results in prime fields. \cite{HZ} obtained a quantitative three-points vector polynomial theorem in $\mathbb R^2$ with two real parameters.
These results make it important to distinguish the many-point, one-parameter, prime field setting of the present theorem.

In one dimension, Bourgain and Chang~\cite{BC} obtained power-saving estimates for the configuration $x,x+t,x+t^2$ in $\F_p$. This was generalized to $x, x+P_1(t), x+P_2(t)$ for arbitrary pairs of linearly independent scalar polynomials $P_1,P_2$ with zero constant terms by Peluse \cite{PeluseThree} and by Dong, Li and Sawin~\cite{DLS}. Peluse~\cite{Peluse} subsequently proved a power-saving counting formula for $x, x+P_1(t),\cdots,x+P_m(t)$ in $\F_p$, when $P_1,\cdots,P_m$ are fixed linearly independent scalar polynomials:
$$
\Lambda_P(f_0,\cdots,f_m) = \prod_{i=0}^m \E f_i + O(p^{-c}).
$$
Thus Corollary~\ref{cor:full-rank} generalizes the above conclusion to vector polynomials whose scalar coordinates are all linearly independent. 

In Peluse's scalar result, the polynomials are fixed before the decay rate; the constants in the present theorem are uniform over all bounded-degree coefficient choices in $\mathbb F_p$.

For multidimensional patterns over $\F_p$, Han, Lacey and Yang~\cite{HLY} obtained power-saving estimates for polynomial corners $x, x+P_1(t)e_1, x+P_2(t)e_2$ (where $x \in \F_p^2$, $t \in \F_p$, and $e_1,e_2$ are the standard basis vectors in $\F_p^2$) when the two polynomials $P_1,P_2$ have distinct degrees, and also when $P_1,P_2$ are linearly independent quadratic polynomials. Kuca~\cite{KucaDistinct} proved a counting formula for configurations of the form
\[
 x,\quad x+P_1(t)v_1,\quad\ldots,\quad x+P_m(t)v_m
\]
with distinct-degree scalar polynomials $P_1,\cdots,P_m$ and nonzero vectors $v_1,\cdots,v_m$. He then replaced the distinct-degree assumption by linear independence of $P_1,\ldots,P_m$ in~\cite{Kuca}. For this type of configurations, each coefficient space $V_i$ is the line $\langle v_i\rangle$ and the independence condition \eqref{eq:independence} becomes the linear independence of $P_1,\cdots,P_m$. Thus Theorem~\ref{thm:main} generalizes Kuca's result to higher dimensional coefficient spaces $V_i$. See also \cite{Lim,LSTW,Hong,HongSim} for further results for corners and for weighted averages.

There are also closely related results on multidimensional polynomial progressions over the integers. 
For example, Peluse, Prendiville, and Shao \cite{PPS} proved an effectir popular counting result for the integer corner $(x,y), (x+t,y), (x,y+t^2)$;
Kravitz, Kuca and Leng~\cite{KKL} develop quantitative concatenation and box norm control for vector-polynomial progressions over $\Z^k$;
Guo, Miao, and Zhan \cite{GMZ} extended distinct-monomial corner estimates to higher-dimensional integer boxes. 
See also \cite{DFKMS,KKS,KMPWW} for related results on vector-polynomial averages in ergodic theory.

\subsection{Proof sketch of Theorem \ref{thm:main}}

Before giving a more detailed proof outline in Section \ref{sec:outline}, we sketch here the key steps in the proof. The proof follows the general strategy of first obtaining higher-order uniformity control and then lowering the order of the norm, developed by Peluse and Prendiville \cite{PP} and by Peluse~\cite{Peluse,PeluseBounds}, and adapted to multidimensional patterns by Kuca~\cite{Kuca}. Our eventual goal is to control the counting operator $\Lambda_P(f_0,\cdots,f_m)$ by the $U^1(V_i)$-norm of $f_i$ for each $1 \leq i \leq m$, where $V_i$ is defind as in \eqref{eq:VW}.
Compared with \cite{Kuca}, the main new difficulty here is that the coefficient spaces $V_i$ need not be lines.

The first step is the box norm estimate in Proposition~\ref{prop:twisted-pet}. Its underlying tools are PET induction, introduced by Bergelson~\cite{Bergelson}, and concatenation~\cite{TZ}, with quantitative developments in~\cite{PP,PeluseBounds}. We use the finite-field estimate of Kuca~\cite[Proposition~4.3]{Kuca}, which already applies to general vector polynomials in our setting and gives control along the lines spanned by the leading coefficient vectors of $P_i-P_j$. We introduce a refinement of this estimate which retains only those lines coming from ``unrestricted'' indices, which is necessary in our inductive argument. This initial uniformity estimate may involve lines which do not lie in the space $V_i$.

We next replace the unwanted lines appearing in the box norm by $V_i$. This is an adaptation of norm smoothing, developed by Frantzikinakis and Kuca~\cite{FK} and used in the finite-field setting in~\cite[Section~5]{Kuca}. This leads to Proposition~\ref{prop:high-norm} which controls the counting operator by the $U^s(V_i)$-norm of each $f_i$ for some positive integer $s$. The new input to this norm smoothing procedure is Proposition \ref{prop:invariant}, a counting result with an invariant factor. For example, consider first the easiest invariant case when $f_1$ is invariant under translation by $V_1$. Then $f_1(x+P_1(t)) = f_1(x)$ for all $x$ and $t$, and hence $\Lambda_P(f_0,\cdots,f_m)$ can be viewed as a counting operator corresponding to the $m-1$ polynomials $P_2,\cdots,P_m$. Thus the desired counting result follows with an appropriate induction setup. In a similar but more complicated manner, Proposition \ref{prop:invariant} proves the desired counting result (with an appropriate induction setup), assuming that one of the functions, say $f_1$, is invariant under translation by some line $L$ (satisfying suitable assumptions). 
The norm-smoothing mechanism  is inherited from \cite{FK}; the additional ingredient here is a counting induction compatible with higher dimensional projections.

Finally, we adapt the degree-lowering argument of Peluse~\cite{Peluse} and Kuca~\cite{Kuca} to pass from control by $U^s(V_i)$-norm to the desired control by $U^1(V_i)$-norm of $f_i$. This is a natural generalization of the arguments in \cite[Section 6]{Kuca} from $1$-dimensional $V_i$'s to general subspaces.

\subsection*{Notation}\label{subsec:averages}

Write $G=\F_p^k$. If $X$ is a finite nonempty set, then $\E_{x\in X}$ denotes the normalized average $|X|^{-1}\sum_{x\in X}$. For a function $f: G\rightarrow\C$ and $1\le q<\infty$, write
\[
 \|f\|_q=(\E_{x\in G}|f(x)|^q)^{1/q},
 \qquad \|f\|_\infty=\max_{x\in G}|f(x)|.
\]
A function is $1$-bounded if its absolute value is at most one. We allow all implied constants to depend on $k,m,d$ (as in the statement of Theorem \ref{thm:main}). Write $e_p(u)=\exp(2\pi i u/p)$ for $u\in\F_p$.

\vspace{0.2in}

\noindent{\bf Use of AI.}
The authors devised the central inductive approach and led the proof's development. AI contributed auxiliary proofs, literature searches and an initial complexity concept that the authors adapted in redesigning the induction. The authors checked every detail and rewrote the manuscript with limited AI-suggested wording, taking full responsibility.

\section{Proof outline}\label{sec:outline}

In this section we give a detailed outline of the proof of Theorem \ref{thm:main}. As mentioned in the introduction,
the proof follows the general strategy of first obtaining a higher-order Gowers norm control for the counting operator and then lowering its order to $U^1$-norms, developed in the works of Peluse~\cite{Peluse,PeluseBounds} and Kuca~\cite{Kuca}. 

\subsection{Counting operators and induction setup}\label{subsec:inductive-theorem}
For $P=(P_1,\ldots,P_m)$ where each $P_i \in \F_p[t]^k$, and functions $f_0,\cdots,f_m: G\rightarrow\C$, recall the definition of the counting operator
$$
\Lambda_P(f_0,\cdots,f_m) = \E_{x \in G} \E_{t \in \F_p} f_0(x) \prod_{i=1}^m f_i(x+P_i(t)).
$$
Thus our main theorem is the statement that $\Lambda_P$ is asymptotically equal to the following simpler linear counting operator obtained by averaging over each shift independently:
\begin{equation}
 M_P(f_0,\ldots,f_m)
 =\E_{x\in G}\E_{v_1\in V_1,\ldots,v_m\in V_m}
       f_0(x)\prod_{j=1}^m f_j(x+v_j) =\E_{x \in G} f_0(x)\prod_{j=1}^m\E(f_j\vert V_j)(x).
\end{equation}

We will prove Theorem \ref{thm:main} by induction, when some of the functions from $f_1,\cdots,f_m$ are arbitrary (unrestricted) and some of them are restricted to be \emph{eigenfunctions}, defined below. This use of eigenfunctions is based on Kuca's degree-lowering argument~\cite[Section~6]{Kuca}. 

\begin{definition}\label{def:eigenfunctions}
For a subspace $V\subset G$, let $\cE(V)$ be the collection of $1$-bounded functions whose restriction to each coset of $V$ is a scalar multiple of a character of $V$. Equivalently, after choosing a complement $Z$ to $V$ so that $G = V\oplus Z$, such a function $\chi$ has the form
\[
 \chi(z+v)=a_ze_p(\lambda_z(v))
\]
for all $z \in Z$ and $v \in V$, where $|a_z| \leq 1$ and $\lambda_z \in V^*$. In this case, we say that $\chi\in \cE(V)$ is an eigenfunction with frequency $\lambda$.
\end{definition}

For a subset $I \subset\{1,\cdots,m\}$, we will consider the counting operator $\Lambda_P(f_0,\cdots,f_m)$, where $f_i$ is an arbitrary $1$-bounded function for $i \in I$, and $f_i \in \cE(V_i)$ for $i \notin I$. For this reason, we refer to $i$ as an \emph{unrestricted index} if $i \in I$ and a \emph{restricted index} if $i \notin I$. Note that if $i \notin I$ is a restricted index and $f_i \in \cE(V_i)$ is an eigenfunction with frequency $\lambda$, then for all $x \in G$ and $v \in V_i$ we have
\begin{equation}\label{eq:eigenfn-property}
f_i(x+v) = f_i(x) e_p(\lambda_x(v)).
\end{equation}
Note also that $\lambda_x \in V_i^*$ depends only on the coset $x+V_i$ in $G/V_i$.

\begin{theorem}\label{thm:strengthened}
Let $P = (P_1,\cdots,P_m)$ be vector polynomials in $\F_p[t]^k$ of degree at most $d$ as in the assumption of Theorem \ref{thm:main}, and let $I \subset \{1,\cdots,m\}$ be a subset. Then for all $1$-bounded functions $f_0,\cdots,f_m:G\rightarrow\C$ with $f_j \in \cE(V_j)$ for $j \notin I$, we have
$$
 \Lam_P(f_0,\cdots,f_m) = M_P(f_0,\cdots,f_m) + O_{k,m,d}(p^{-c})
$$
for some constant $c = c(k,m,d) > 0$.
\end{theorem}

Clearly Theorem \ref{thm:main} is a special case of Theorem \ref{thm:strengthened} with $I = \{1,\cdots,m\}$. 
We will prove Theorem \ref{thm:strengthened} by induction, using the following notion of complexity of $(P,I)$.

\begin{definition}[Complexity of $(P,I)$]\label{def:complexity}
Let $P = (P_1,\cdots,P_m)$ be vector polynomials in $\F_p[t]^k$ with $P_i(0) = 0$ for each $1 \leq i \leq m$, and let $I \subset \{1,\cdots,m\}$ be a subset. Recall the definition of $V_i$ from \eqref{eq:VW}. Define the complexity of $(P,I)$ to be the quantity
$$
\mathfrak{c}(P,I) = m + |I| + 2\sum_{i \in I}\dim V_i + \dim\Big(\sum_{i \in I}V_i\Big).
$$
\end{definition}

Thus the complexity of $(P,I)$ depends on $m$, on the number of unrestricted index, and on the dimensions of the spaces $V_i$ for those $i \in I$. In our induction process, we will frequently convert a counting operator $\Lambda_P$ associated with $(P,I)$ to a new counting operator via the following processes, all of which strictly reduce the complexity:
\begin{enumerate}
\item Removing a polynomial from $P$;
\item Moving one unrestricted index $i \in I$ to restricted;
\item Projecting each $P_i$ for $i \in I$ to $TP_i$, where $T$ is a projection from $\sum_{i\in I}V_i$ to a proper subspace.
\item Projecting one polynomial $P_i$ for some $i \in I$ to $TP_i$, where $T$ is a projection from $V_i$ to a proper subspace, while creating a new polynomial with a restricted index;
\end{enumerate}

\subsection{Initial control by higher-order norms}

The first step in the proof of Theorem \ref{thm:strengthened} is to obtain control for our counting operator by a higher-order box norm using PET induction~\cite{Bergelson} and concatenation~\cite{TZ}. Quantitative versions of these processes (in the integers) were initially developed by Peluse and Prendiville~\cite{PP} and in Peluse's work~\cite{PeluseBounds}.

\begin{proposition}[Box norm control]\label{prop:twisted-pet}
Let $P = (P_1,\cdots,P_m)$ be  nonzero vector polynomials in $\F_p[t]^k$ of degree at most $d$ as in the assumption of Theorem \ref{thm:main}, let $I \subset \{1,\cdots,m\}$ be a subset, and let $i \in I$. For $j \in (I\setminus\{i\}) \cup \{0\}$, let
$$
L_{ij} = \langle \lc(P_i-P_j)\rangle \subset G,
$$
where $\langle v\rangle$ denotes the span of a vector $v \in G$, $\lc(P_i-P_j)$ denotes the largest coefficient of $P_i-P_j$, and we set $P_0 = 0$. Then for all $1$-bounded functions $f_0,\cdots,f_m:G\rightarrow\C$ with $f_j \in \cE(V_j)$ for $j \notin I$, we have
$$
 |\Lam_P(f_0,\cdots,f_m)|\ll_{k,m,d} \|f_i\|_{\mathcal B_i}^{c}+p^{-c}
$$
for some constant $c = c(k,m,d) > 0$, where $\mathcal B_i$ is a (repeated) list of $O_{k,m,d}(1)$ lines from $\{L_{ij}\}$ with $j \in (I\setminus\{i\})\cup\{0\}$.
\end{proposition}

The definition of these box norms are given in Section \ref{sec:prelim}.
Unfortunately Proposition \ref{prop:twisted-pet} does not give control along the ``correct'' space $V_i$ as needed in our main theorem. For example, take
\begin{equation}\label{eq:outline-example}
 P_1(t)=e_1t+e_2t^2,\qquad P_2(t)=e_3t^3,
\end{equation}
where $e_1,e_2,e_3$ are the standard basis vectors in $\F_p^3$.
Then the assumption \eqref{eq:independence} is satisfied. On the other hand, we have
\[
 V_1=\langle e_1,e_2\rangle,
 \qquad \lc(P_1-P_2)=-e_3.
\]
Thus Proposition \ref{prop:twisted-pet} gives control for a box norm of $f_1$ along $e_3$, whereas we eventually want to control $f_1$ along $V_1$.

Proposition \ref{prop:twisted-pet} will be proved in Section \ref{sec:PET}. 
Its proof starts from Kuca's box norm estimate for vector polynomials, stated below as Proposition \ref{prop:pet-box}. The additional issue is to ensure that the lines $L_{ij}$ appearing in the box norm only involve unrestricted indices $j$ or $j=0$. This is achieved by removing the eigenfunctions at restricted indices via repeated Cauchy--Schwarz, before applying Kuca's box norm estimate.

\subsection{Higher-order control along $V_i$ via norm smoothing}

Next we replace the lines in $\mathcal B_i$ by the coefficient space $V_i$, using the norm smoothing argument of \cite{FK,Kuca}.

\begin{proposition}[Higher-order control along $V_i$]\label{prop:high-norm}
Let $P = (P_1,\cdots,P_m)$ be vector polynomials in $\F_p[t]^k$ of degree at most $d$ as in the assumption of Theorem \ref{thm:main}, and let $I \subset \{1,\cdots,m\}$ be a  nonempty subset.
 Suppose that Theorem \ref{thm:strengthened} has already been established for $(P',I')$ with $\mathfrak c(P',I') < \mathfrak c(P,I)$.
Then for all $1$-bounded functions $f_0,\cdots,f_m:G\rightarrow\C$ with $f_j \in \cE(V_j)$ for $j \notin I$, we have
$$
 |\Lam_P(f_0,\cdots,f_m)|\ll_{k,m,d} \min_{i\in I}\|f_i\|_{U^s(V_i)}^{c}+p^{-c}
$$
for some constant $c = c(k,m,d) > 0$ and positive integer $s = s(k,m,d)$.
\end{proposition}

The additional input  in the proof of Proposition \ref{prop:high-norm} is a counting result when one unrestricted function is invariant along a line.

\begin{proposition}[Counting in the $L$-invariant case]\label{prop:invariant}
Let $P = (P_1,\cdots,P_m)$ be vector polynomials in $\F_p[t]^k$ of degree at most $d$ as in the assumption of Theorem \ref{thm:main}, and let $I \subset \{1,\cdots,m\}$ be a subset. Suppose that Theorem \ref{thm:strengthened} has already been established for $(P',I')$ with $\mathfrak c(P',I') < \mathfrak c(P,I)$. Let $L \subset \sum_{i\in I}V_i$ be a line. Then for all $1$-bounded functions $f_0,\cdots,f_m: G\rightarrow\C$ with $f_j \in \cE(V_j)$ for $j \notin I$ and with $f_{i}$ $L$-invariant for some $i \in I$, we have
$$
\Lambda_P(f_0,\cdots,f_m) = M_P(f_0,\cdots,f_m) + O_{k,m,d}(p^{-c})
$$
for some constant $c=c(k,m,d) > 0$.
\end{proposition}

To illustrate the argument of deducing Proposition \ref{prop:high-norm} from Propositions \ref{prop:twisted-pet} and \ref{prop:invariant}, consider the example \eqref{eq:outline-example} with $I=\{1,2\}$. Here $V_1=\langle e_1,e_2\rangle$ and $V_2=\langle e_3\rangle$, while $\mathcal B_1$ consists of repeated copies of $L_{01}=\langle e_2\rangle$ and $L_{21}=\langle e_3\rangle$. First write $\mathcal B_1=(\mathcal B,\langle e_3\rangle)$ and we aim at replacing $\langle e_3\rangle$ by $V_1$. Passing to the dual function associated with $f_1$ and applying dual-difference interchange (Lemma \ref{lem:dual-difference}) produces counting operators of the form $\Lambda_P(g_0,g_1,g_2)$ where $g_1$ is $\langle e_3\rangle$-invariant. Proposition \ref{prop:invariant} applies to this operator, giving
$$
|\Lambda_P(g_0,g_1,g_2)| \leq \E_{x \in G} |\E(g_2\vert V_2)(x)|+O(p^{-c}).
$$
By using the specific shape of the function $g_2$, this eventually leads to 
\[
 |\Lambda_P(f_0,f_1,f_2)|
 \ll \|f_2\|_{\mathcal B,V_2}^{c}+p^{-c}.
\]
This process of transferring from control by the box norm $\|f_1\|_{\cB,\langle e_3\rangle}$ to control by the box norm $\|f_2\|_{\cB,V_2}$ is summarized in Lemma \ref{lem:one-transfer}.
We then repeat the transfer process above from control by $\|f_2\|_{\cB,V_2}$ to control by $\|f_1\|_{\cB,V_1}$.
This replaces one line in $\mathcal B_1$ by $V_1$. Repeating the two transfers replaces every line in $\mathcal B_1$ by $V_1$, giving $U^s(V_1)$ control. The details are given in Section \ref{sec:smoothing}.

To illustrate the argument of proving Proposition \ref{prop:invariant}, consider again the previous example \eqref{eq:outline-example} with $I = \{1,2\}$. Suppose that $f_1$ is $L$-invariant with $L=\langle e_3\rangle$, a situation that arises in the discussion above. By shifting $x$ by $\ell$ and using the $L$-invariance of $f_1$, we have
$$
\Lambda_P(f_0,f_1,f_2) = \E_{x,t} f_1(x+P_1(t)) \E_{\ell \in L} f_0(x+\ell)f_2(x+
\ell+P_2(t)).
$$
By Cauchy-Schwarz in the $\ell$ variable, we have
$$
|\Lambda_P(f_0,f_1,f_2)|^2 \leq \E_{\ell,\ell'\in L} \E_{x,t} f_0(x+\ell)\overline{f_0(x+\ell')} f_2(x+\ell+P_2(t))\overline{f_2(x+\ell'+P_2(t))}+O(p^{-c}).
$$
For fixed $\ell,\ell'$, the inner average over $x,t$ can be controlled using the $m=1$ case of Theorem \ref{thm:strengthened} (or by induction hypothesis in the general case). This gives
$$
|\Lambda_P(f_0,f_1,f_2)|^2  \leq \E_{\ell,\ell'\in L} \E_{x \in G} \E_{v \in V_2} f_0(x+\ell)\overline{f_0(x+\ell')} f_2(x+\ell+v)\overline{f_2(x+\ell'+v)}+O(p^{-c}).
$$
By Cauchy-Schwarz in the $v$ variable to remove $f_0$, the square of the right-hand side above is at most the box norm $\|f_2\|_{V_2,L}$. This control by a mixed box norm is summarized in Lemma \ref{lem:invariant-energy}.
On the other hand, the elementary inverse theorem for the box norm $\|\cdot\|_{V_2,L}$ (Lemma \ref{lem:invariant-mixed-inverse}), applied to the appropriate dual function, allows us to replace $f_2$ in the counting operator by a structured function $K$. In this example we have $L=V_2$, and $K$ is a product of a $L$-invariant function $h$ and an eigenfunction $\chi \in \cE(V_2)$. Then
$$
\begin{aligned}
\Lambda_P(f_0,f_1,K) &= \E_{x,t} f_0(x) f_1(x+P_1(t)) h(x+P_2(t)) \chi(x+P_2(t))\\
&= \E_{x,t} (f_0h)(x) f_1(x+P_1(t)) \chi(x+P_2(t)) = \Lambda_P(f_0h, f_1,\chi).
\end{aligned}
$$
The latter counting operator corresponds to $(P, \{1\})$ which has smaller complexity, and hence the induction hypothesis applies. The details are given in Section \ref{sec:invariant}.

\subsection{Degree lowering: from $U^s(V_i)$ to $U^1(V_i)$}

The final step is to pass from control by $U^s(V_i)$-norm to control by $U^1(V_i)$-norm, following the degree-lowering argument of \cite{Peluse,Kuca}. To illustrate the argument of deducing Theorem \ref{thm:strengthened} from Proposition \ref{prop:high-norm}, again we use example \eqref{eq:outline-example} with $I = \{1,2\}$. Proposition \ref{prop:high-norm} already gives
$$
|\Lambda_P(f_0,f_1,f_2)| \ll \|f_1\|_{U^s(V_1)}^c + p^{-c}
$$
for some positive integer $s$. If $s \geq 3$, we perform the following procedure (Lemma \ref{lem:degree-lowering}) which effectively replaces $s$ by $s-1$ in the estimate above. Passing to the dual function and applying dual-difference interchange (Lemma \ref{lem:dual-difference}) produces counting operators of the form $\Lambda_P(g_0,g_1,g_2)$ where $g_1 \in \cE(V_1)$ is an eigenfunction. The induction hypothesis applies to this operator, giving
$$
|\Lambda_P(g_0,g_1,g_2)| \leq \E_{x \in G} |\E(g_1\vert V_1)(x)|+O(p^{-c}).
$$
By using the formula for the function $g_1$, this eventually leads to
$$
|\Lambda_P(f_0,f_1,f_2)| \ll \|f_1\|_{U^{s-1}(V_1)}^c + p^{-c}.
$$
Repeating the process above eventually lowers the order to $U^2(V_1)$, at which point a standard Fourier analysis argument can further lower the order to $U^1(V_1)$. The counting formula in Theorem \ref{thm:strengthened} then follows from a standard telescoping argument. The details are given in Section \ref{sec:degree-lowering}.

We end this section with a diagram illustrating the various components of the proof and where they are proven in the paper. Solid arrows show the flow of the argument, and dashed arrows indicate uses of the induction hypothesis in Theorem \ref{thm:strengthened} (counting theorem at smaller complexity).

\begin{center}
\begin{minipage}{\textwidth}
\centering
\begin{tikzpicture}[
  >=Stealth,
  result/.style={draw,rounded corners=2pt,align=center,text width=4.5cm,
    minimum height=1.05cm,inner sep=5pt,font=\small},
  flow/.style={->,semithick},
  induction/.style={->,semithick,dashed}]
\node[result] (pet) at (-3.25,0)
  {Proposition~\ref{prop:twisted-pet} (Section \ref{sec:PET})\\Box norm control\\PET and concatenation};
\node[result] (hyp) at (3.25,0)
  {Induction hypothesis\\ in Theorem~\ref{thm:strengthened}};
\node[result] (smooth) at (-3.25,-2.15)
  {Proposition~\ref{prop:high-norm} (Section \ref{sec:smoothing})\\$U^s(V_i)$ control\\Norm smoothing};
\node[result] (invariant) at (3.25,-2.15)
  {Proposition~\ref{prop:invariant} (Section \ref{sec:invariant})\\Counting with an\\$L$-invariant function};
\node[result] (lower) at (3.25,-4.4)
  {Lemma~\ref{lem:degree-lowering} (Section \ref{sec:degree-lowering})\\ Degree lowering \\$U^s(V_i)\longrightarrow U^2(V_i)$};
\node[result,text width=6.0cm] (count) at (0,-7.0)
  {Theorem~\ref{thm:strengthened}:   $\Lambda_P=M_P+O(p^{-c})$};
\node[result,text width=4.5cm,minimum height=.8cm] (main) at (0,-8.6)
  {Theorem~\ref{thm:main}};
\draw[flow] (pet) -- (smooth);
\draw[induction] (hyp) -- (invariant);
\draw[induction] (hyp.south west) -- (smooth.north east);
\draw[flow] (invariant) -- (smooth);
\draw[flow] (smooth.south) -- (-3.25,-4.4) -- (count);
\draw[induction] (hyp.east) -- (6.0,0) -- (6.0,-4.4) -- (lower.east);
\draw[induction] (6.0,-4.4) -- (6.0,-7.0) --(count.east);
\draw[flow] (lower.south) -- (3.25,-5.65) -| (count.north);
\draw[flow] (count) -- node[right,font=\footnotesize]{$I=\{1,\ldots,m\}$} (main);
\end{tikzpicture}
\end{minipage}
\end{center}

\section{Preliminaries}\label{sec:prelim}

\subsection{Box norms and Fourier analysis}
Recall that $G = \F_p^k$.
For a function $f: G\rightarrow \C$ and $h\in G$, write $\Delta_hf(x)=f(x)\overline{f(x+h)}$. For a tuple $h=(h_1,\cdots,h_s)$, write
\[
 \Delta_hf(x)=\Delta_{h_1}\cdots\Delta_{h_s}f(x)
 =\prod_{\omega\in\{0,1\}^s}\cC^{|\omega|}
       f(x+\omega_1h_1+\cdots+\omega_sh_s),
\]
where $\cC$ denotes complex conjugation.
For a list $\mathcal B=(H_1,\ldots,H_s)$ of subgroups of $G$ with $s\geq1$, define the box norm $\|f\|_{\cB}$ by
\[
\begin{aligned}
 \|f\|_{\mathcal B}^{2^s} &= \E_{x \in G} \E_{h \in H_1\times\cdots\times H_s} \Delta_hf(x) \\
 &=\E_{x\in G}\E_{h_1\in H_1,\ldots,h_s\in H_s}
   \prod_{\omega\in\{0,1\}^s}\cC^{|\omega|}
       f(x+\omega_1h_1+\cdots+\omega_sh_s).
\end{aligned}
\]
See \cite[Section~2.1]{Kuca} for some basic properties of these box norms. For example, if $\cB'$ is obtained by adding entries to $\cB$ then $\|f\|_{\cB} \leq \|f\|_{\cB'}$.

For a subspace $V\subset G$, define the $U^s(V)$-norm by  $\|f\|_{U^s(V)}=\|f\|_{V^{\times s}}$, where $V^{\times s} = (V, \cdots, V)$ denotes a list of $s$ copies of $V$. In particular,
\begin{equation}
 \|f\|_{U^1(V)}=\|\E(f\vert V)\|_2.
\end{equation}
If $\mathcal B=(H_1,\cdots,H_s)$ and $V\subset G$ is another subspace, expanding the definition gives
\begin{equation}
\begin{aligned}
 \|f\|_{\mathcal B,V}^{2^{s+1}} &= \E_{h \in H_1\times\cdots\times H_s} \E_{x \in G} \E_{v \in V} \Delta_v\Delta_hf(x) \\
&= \E_{h \in H_1\times \cdots\times H_s} \|\Delta_hf\|_{U^1(V)}^2
 =\E_{h\in H_1\times\cdots\times H_s}
       \|\E(\Delta_hf\vert V)\|_2^2.
\end{aligned}
\end{equation}

For a subspace $V \subset G$ and  a function $g:V\rightarrow\C$, define its normalized Fourier transform $\wh g: V^*\rightarrow \C$ by
\[
 \wh g(\xi)=\E_{v\in V}g(v)e_p(-\xi(v))
\]
for $\xi \in V^*$.
Then Fourier inversion and Parseval's identity give
$$
g(v)=\sum_{\xi\in V^*}\wh g(\xi)e_p(\xi(v)), \ \
\sum_{\xi\in V^*}|\wh g(\xi)|^2=\E_{v\in V}|g(v)|^2.
$$

Choose a complement $Z$ with $G=V\oplus Z$. For a function $f: G\rightarrow\C$, write $f_z(v)=f(z+v)$ for $z\in Z$ and $v\in V$. We have
\begin{equation}
 \|f\|_{U^s(V)}^{2^s}
 =\E_{z\in Z}\|f_z\|_{U^s(V)}^{2^s},
 \qquad
 \|f\|_{U^2(V)}^4
 =\E_{z\in Z}\sum_{\xi\in V^*}|\wh{f_z}(\xi)|^4.
\end{equation}
More generally, for $s\geq2$,
\begin{equation}
 \|f\|_{U^s(V)}^{2^s}
 =\E_{z\in Z}\E_{h\in V^{s-2}}
       \sum_{\xi\in V^*}|\wh{\Delta_hf_z}(\xi)|^4.
\end{equation}
If $f$ is $1$-bounded, then these identities and Parseval imply
\begin{equation}
 \|f\|_{U^2(V)}^4
 \leq\E_{z\in Z}\max_{\xi\in V^*}|\wh{f_z}(\xi)|
 =\sup_{\chi\in\cE(V)}\left|\E_{x\in G}f(x)\chi(x)\right|,
\end{equation}
where $\chi \in \cE(V)$ is defined by $\chi(z+v)=a_ze_p(-\xi_z(v))$ for $z \in Z$ and $v \in V$, with $|a_z|=1$ and $\xi_z \in V^*$ chosen such that $a_z\wh{f_z}(\xi_z)=\max_{\xi\in V^*}|\wh{f_z}(\xi)|$.

For $s\geq2$ and $g: \F_p\rightarrow\C$, we also use the dual norm
\[
 \|g\|_{U^s(\F_p)^*}
 =\sup_{f}|\E_{t\in\F_p}f(t)g(t)|,
\]
where the supremum is over all functions $f:\F_p\rightarrow\C$ satisfying $\|f\|_{U^s(\F_p)} \leq 1$.
If $P\in\F_p[t]$ is a polynomial of degree less than $s$, then
$\Delta_h e_p(P)=1$ for every $h\in\F_p^s$, 
which implies also that $\|e_p(P(t))\|_{U^s(\mathbb F_p)}=1$. Hence multiplication by $e_p(P)$ preserves the $U^s(\F_p)$-norm, and
\[
 |\E_t f(t)e_p(P(t))|\leq\|fe_p(P)\|_{U^s(\F_p)}
 =\|f\|_{U^s(\F_p)}.
\]
Conversely, taking $f(t)=e_p(-P(t))$ implies that $\|e_p(P)\|_{U^s(\F_p)^*}\ge1$.
In summary $\|e_p(P)\|_{U^s(\F_p)^*}=1$.

\subsection{Projection onto a prescribed frequency}
Recall the definition of $\cE(V)$ from Definition \ref{def:eigenfunctions}.
The class $\cE(V)$ is closed under translations, complex conjugation and products. In particular, $\Delta_hf\in\cE(V)$ whenever $f\in\cE(V)$ and $h\in G^s$.

\begin{lemma}\label{lem:frequency-projection}
Let $V\subset G$ be a subspace, and for each $x \in G$ let $\lambda_x \in V^*$ be a frequency depending only on the coset $x+V$ in $G/V$. For a $1$-bounded function $f: G\rightarrow\C$, define $\Pi_\lambda f: G\rightarrow \C$ by
\[
 \Pi_\lambda f(x)=\E_{v\in V}f(x+v)e_p(-\lambda_x(v))
\]
for $x \in G$. 
\begin{enumerate}
\item We have $\Pi_\lambda f\in\cE(V)$ with frequencey $\lambda$ and $\langle f, \Pi_{\lambda}f\rangle = \|\Pi_{\lambda}f\|_2^2$.
\item For any  $\chi \in \cE(V)$ with frequency $\lambda$, we have $\langle \chi,f\rangle = \langle \chi, \Pi_{\lambda}f\rangle$.
\item For any list $\mathcal B$ of subspaces, we have $\|\Pi_\lambda f\|_{\mathcal B, V}\leq\|f\|_{\mathcal B, V}$.
\end{enumerate}
\end{lemma}

\begin{proof}
Writing $G=V\oplus Z$, we have
$$
\Pi_\lambda f(z+w) = \E_{v \in V} f_z(w+v) e_p(\lambda_z(v)) = \wh{f_z}(\lambda_z) e_p(\lambda_z(w))
$$ 
for $z \in Z$ and $w \in V$. This shows that $\Pi_\lambda \in cE(V)$ has frequency $\lambda$. Moreover, we have
$$
\langle f,\Pi_\lambda f\rangle = \E_{z \in Z, w\in V} f_z(w) \overline{\wh{f_z}(\lambda_z)} e_p(-\lambda_z(w)) = \E_{z \in Z} |\wh{f_z}(\lambda_z)|^2
$$
and
$$
\|\Pi_{\lambda}f\|_2^2 = \E_{z \in Z, w\in V} |\Pi_{\lambda}f(z+w)|^2 = \E_{z \in Z} |\wh{f_z}(\lambda_z)|^2.
$$
This proves (1). For (2), writing $\chi(z+w) = a_ze_p(\lambda_z(w))$ for some $a_z \in \C$ with $|a_z| \leq 1$, we have
$$
\langle \chi,\Pi_\lambda f\rangle = \E_{z \in Z} a_z\overline{\wh{f_z}(\lambda_z)} = \E_{z \in Z} a_z \E_{v \in V} \overline{f_z(v)} e_p(\lambda_z(v)) = \langle\chi,f\rangle.
$$
For (3), by the triangle inequality for the box norm we have
$$
\|\Pi_{\lambda}f\|_{\cB,V} \leq \E_{v \in V} \|f(\cdot + v) g_v\|_{\cB,V},
$$
where $g_v: G\rightarrow \C$ denotes the function $g_v(x) = e_p(-\lambda_x(v))$. Then $g_v$ is $V$-invariant: $g_v(x+w) = g_v(x)$ for any $w \in V$. By the definition of the box norm, if $\cB = (H_1,\cdots,H_{s-1})$ then
$$
\|f(\cdot + v) g_v\|_{\cB,V}^{2^{s}} = \E_{x \in G} \E_{h \in H_1\times\cdots\times H_{s-1}} \E_{w \in V} \Delta_h\Delta_w f(x+v) \Delta_h\Delta_w g_v(x).
$$
Since $\Delta_wg_v(x) = g_v(x)\overline{g_v(x+w)} = |g_v(x)|^2 =1$ by the $V$-invariance of $g_v$, we have
$$
\|f(\cdot + v) g_v\|_{\cB,V}^{2^{s}} = \E_{x \in G} \E_{h \in H_1\times\cdots\times H_{s-1}} \E_{w \in V} \Delta_h\Delta_w f(x+v) = \|f(\cdot +v)\|_{\cB,V}^{2^{s}}.
$$
The claim follows since the box norm is translation invariant.
\end{proof}

\subsection{Boundary cases of Theorem \ref{thm:strengthened}}

We use the classical Weil bound \cite{Weil}: if $P(t)\in\F_p[t]$ is nonconstant polynomial of degree $d<p$, then
\begin{equation}\label{eq:prelim-Weil}
 \left|\E_{t\in\F_p}e_p(P(t))\right|\leq(d-1)p^{-1/2}.
\end{equation}

\begin{lemma}\label{lem:all-restricted}
Let $P=(P_1,\cdots,P_m)$ be vector polynomials in $\F_p[t]^k$ of degree at most $d$ as in the assumption of Theorem \ref{thm:main}. For $1$-bounded functions $f_0,\cdots,f_m:G\rightarrow\C$ with $f_j\in\cE(V_j)$ for every $1\leq j\leq m$, we have
\[
 \Lambda_P(f_0,\cdots,f_m)=M_P(f_0,\cdots,f_m)+O_d(p^{-1/2}).
\]
\end{lemma}

\begin{proof}
For $1 \leq j \leq m$, choose frequencies $\lambda_{j,x}\in V_j^*$ as in \eqref{eq:eigenfn-property}. Then
$$
f_j(x+P_j(t)) = f_j(x) e_p(\lambda_{j,x}(P_j(t)))
$$
for all $x \in G$ and $t\in\F_p$. We have
$$
\Lambda_P(f_0,\cdots,f_m) = \E_{x \in G} (f_0f_1\cdots f_m)(x) \E_{t\in\F_p} e_p\Big(\sum_{j=1}^m \lambda_{j,x}(P_j(t))\Big).
$$
The total contribution from those $x \in G$ such that $\sum_j \lambda_{j,x}(P_j(t))$ is nonconstant is $O_d(p^{-1/2})$ by the Weil bound \eqref{eq:prelim-Weil}. 

Now suppose that $\sum_j \lambda_{j,x}(P_j(t))$ is constant. This constant has to be $0$ since $P_j(0) = 0$ for every $j$. Since $\lambda_{j,x}(P_j(t)) \in W_j$ by the definition of $W_j$ in \eqref{eq:VW}, the direct sum assumption \eqref{eq:independence} implies that $\lambda_{j,x}(P_j(t)) = 0$ for every $1 \leq j \leq m$. Since the coefficients of $P_j$ span $V_j$, this implies that $\lambda_{j,x}=0$ for every $1 \leq j\leq m$. It follows from that
$$
\Lambda_P(f_0,\cdots,f_m) = \E_{x \in G} f_0(x) \prod_{1 \leq j \leq m} f_j(x) 1_{\{\lambda_{j,x}=0\}} + O_d(p^{-1/2}).
$$
On the other hand, note that
$$
\E(f_j\vert V_j)(x) = \E_{v \in V_j}f_j(x+v) = f_j(x) \E_{v \in V_j}  e_p(\lambda_{j,x}(v)) = f_j(x) 1_{\{\lambda_{j,x}=0\}}.
$$
This completes proof by recalling the definition of $M_P(f_0,\cdots,f_m)$.
\end{proof}

\begin{lemma}\label{lem:one-shift}
Let $P\in\F_p[t]^k$ be a vector polynomial with degree at most $d$ and $P(0)=0$. Let $V=V(P)$ as in \eqref{eq:VW}. For any function $f:G\rightarrow\C$, we have
$$
 \left\|\E_{t\in\F_p}f(\cdot+P(t))-\E(f\vert V)\right\|_2
 \ll_d p^{-1/2}\|f\|_2.
$$
\end{lemma}

\begin{proof}
The left-hand side is the $L^2$-norm of the function $F(x) = \E_t f(x+P(t)) - \E_{v \in V} f(x+v)$, whose Fourier transform can be computed to be
$$
\begin{aligned}
\wh{F}(\xi) &= \E_{x\in G}\E_{t \in \F_p} f(x+P(t)) e_p(-\xi(x)) - \E_{x \in G} \E_{v \in V}f(x+v) e_p(-\xi(x)) \\
&= \E_{x \in G} f(x) e_p(-\xi(x)) \Big(\E_{t \in \F_p}e_p(\xi(P(t))) - \E_{v \in V} e_p(\xi(v))\Big).
\end{aligned}
$$
If $\xi$ is identically zero on $V$, then both the average over $t \in \F_p$ and the average over $v \in V$ above are equal to $1$, and hence the difference vanishes. Otherwise, the average over $v \in V$ is zero and the Weil bound \eqref{eq:prelim-Weil} implies that $\E_t e_p(\xi(P(t)) \ll_d p^{-1/2}$. It follows that
$$
|\wh{F}(\xi)| \ll_d p^{-1/2} |\wh{f}(\xi)|
$$
for all $\xi$. Hence $\|F\|_2^2 = \|\wh{F}\|_2^2 \ll_d p^{-1} \|\wh{f}\|_2^2 = p^{-1}\|f\|_2^2$.
\end{proof}

In particular, Lemma~\ref{lem:one-shift} proves Theorem~\ref{thm:strengthened} when $m=1$: 
$$
\begin{aligned}
|\Lambda_P(f_0,f_1) - M_P(f_0,f_1)| &= \Big|\E_{x \in G}f_0(x) \Big(\E_{t \in \F_p}f_1(x+P_1(t)) - \E(f_1\vert V_1)(x)\Big)\Big| \\
&\leq \left\|\E_{t\in \F_p}f_1(\cdot +P_1(t)) - \E(f_1\vert V_1)\right\|_2 \ll_d p^{-1/2}.
\end{aligned}
$$
In view of Lemmas \ref{lem:all-restricted} and \ref{lem:one-shift}, we can always assume that $I$ is nonempty and $m \geq 2$ when proving the statements in Section \ref{sec:outline}.
In particular, we may assume $d\ge 2$, since by the block independence condition \eqref{eq:independence}, we have $m=1$ when all polynomials are linear.

\subsection{Cauchy--Schwarz lemmas}
We record subgroup versions of the two Cauchy--Schwarz lemmas in \cite[Lemmas~A.2 and~A.4]{Kuca}. 

\begin{lemma}\label{lem:missing-variable-CS}
Let $\mathcal B=(H_1,\cdots,H_s)$ be a list of subgroups of $G$, and let $f:G\rightarrow\C$ be $1$-bounded. For each $1\leq\ell\leq s$ and $h=(h_1,\cdots,h_s)\in H_1\times\cdots\times H_s$, let $g_{\ell,h}: G\rightarrow \C$ be a $1$-bounded function. Suppose that $g_{\ell,h}$ does not depend on $h_\ell$ (i.e. if $h = (h_1,\cdots,h_s)$ and $h'=(h_1',\cdots,h_s')$ satisfy $h_i = h_i'$ for all $i \neq\ell$ then $g_{\ell,h} = g_{\ell,h'}$). Then
$$
 \left|\E_{x \in G}\E_{h \in H_1\times\cdots\times H_s}\Delta_hf(x)\prod_{\ell=1}^s g_{\ell,h}(x)\right|
 \leq\|f\|_{\mathcal B}.
$$
\end{lemma}

\begin{proof}
For $h,h'\in H_1\times\cdots\times H_s$ and $\omega\in\{0,1\}^s$, define $h^\omega = (h_1^\omega,\cdots,h_s^\omega)$ by
\begin{equation}\label{eq:prelim-homega}
h^\omega_\ell = \begin{cases}  h_{\ell} & \text{if }\omega_\ell=0, \\ h'_\ell & \text{if }\omega_\ell=1. \end{cases}
\end{equation}
Since each factor in the product
$$
\Delta_hf(x) = \prod_{\omega \in \{0,1\}^s} \cC^{|\omega|} f(x+\omega_1h_1+\cdots+\omega_sh_s)
$$
(except for $\cC^s f(x+h_1+\cdots+h_s)$) is independent of at least one coordinate of $h$, we can write
$$
\Delta_hf(x) \prod_{\ell=1}^s g_{\ell,h}(x) = b(h) \prod_{\ell=1}^s a_{\ell}(h),
$$
where $b(h) = \cC^s f(x+h_1+\cdots+h_s)$ and each $a_{\ell}$ is $1$-bounded and independent of $h_{\ell}$.
By successively applying the Cauchy--Schwarz inequality to remove $a_1,\cdots,a_s$, we have
\begin{equation}\label{eq:coordinate-CS}
 \left|\E_h b(h)\prod_{\ell=1}^s a_\ell(h)\right|^{2^s}
 \leq\E_{h,h'}\prod_{\omega\in\{0,1\}^s}\cC^{|\omega|}b(h^\omega).
\end{equation}
By H\"{o}lder's inequality, the $2^s$th power of the left-hand side of  the desired inequality is at most
\[
 \E_x\E_{h,h'}\prod_{\omega\in\{0,1\}^s}
       \cC^{s+|\omega|}f(x+h^\omega_1+\cdots+h^\omega_s)
 =\E_x\E_{h,h'}\Delta_{h-h'}f(x)
 =\|f\|_{\mathcal B}^{2^s}.
\]
This completes the proof.
\end{proof}

\begin{lemma}[Dual-difference interchange]\label{lem:dual-difference}
Let $P=(P_1,\cdots,P_m)$ be vector polynomials in $\F_p[t]^k$, let $P_0=0$, and let $1\leq i\leq m$. For $1$-bounded functions $f_j:G\rightarrow\C$ with $0\leq j\leq m$ and $j\ne i$, define
\[
 d_i(x)=\E_{t \in \F_p}\prod_{\substack{0\leq j\leq m\\j\ne i}}
                  f_j(x+P_j(t)-P_i(t)).
\]
Let $H_1,\cdots,H_s\subset G$ be subspaces, and for each $h\in H_1\times\cdots\times H_s$ let $u_h:G\rightarrow\C$ be a $1$-bounded function. For $h,h' \in H_1\times\cdots\times H_s$, let $U_{h,h'}:G\rightarrow\C$ be defined by
\[
 U_{h,h'}(x)=\prod_{\omega\in\{0,1\}^s}\cC^{|\omega|}
                u_{h^\omega}\Bigl(x-\sum_{\ell=1}^s h'_\ell\Bigr),
\]
where $h^\omega$ is defined as in \eqref{eq:prelim-homega}.
Then
$$
 \left|\E_h\E_x\Delta_hd_i(x)u_h(x)\right|^{2^s}
 \leq\E_{h,h'}\left|\Lambda_P
       (\Delta_{h-h'}f_0,\cdots,\Delta_{h-h'}f_{i-1}, U_{h,h'},\Delta_{h-h'}f_{i+1}\cdots,\Delta_{h-h'}f_m)\right|,
$$
where $U_{h,h'}$ occupies the $i$th component.
\end{lemma}

\begin{proof}
Write $g_t(x)=\prod_{j\ne i}f_j(x+P_j(t)-P_i(t))$, so that $d_i=\E_tg_t$ and we have
$$
\E_h\E_x\Delta_hd_i(x)u_h(x) = 
\E_{x} \E_h [\E_t \mathcal C^{s}g_t(x+\sum_{1\le\ell s} h_\ell) u_h(x)][\prod_{\substack{\omega\in\{0,1\}^s,\\\omega\neq(1,\dots,1)}} \E_t \mathcal C^{|\omega|}g_t(x+\omega\cdot h)] .
$$

Since each factor $\E_t \mathcal C^{|\omega|}g_t(x+\omega\cdot h)$ with $\omega\neq(1,\dots,1)$ in the last product
independent of at least one coordinate of $h$, 
we can write the product as
$$
\prod_{\ell=1}^s a_{\ell}(h),
$$
where %$b(h) = u_h(x) \cC^s g_t(x+h_1+\cdots+h_s)$ and 
each $a_{\ell}$ is $1$-bounded and independent of $h_{\ell}$. 
By \eqref{eq:coordinate-CS}, we have
$$
|\E_h \mathcal C^{s}g_t(x+\sum_{1\le\ell s} h_\ell) u_h(x)\prod_{\ell=1}^s a_{\ell}(h) |^{2^s} 
\leq \E_{h,h'} \prod_{\omega \in\{0,1\}^s} \cC^{|\omega|} u_{h^\omega}(x) \cdot\cC^{|\omega|+s} g_t(x+h_1^\omega+\cdots+h_s^\omega).
$$
By H\"{o}lder's inequality, the left-hand side of the desired inequality is at most
$$
\E_{x,t} \E_{h,h'}\prod_{\omega \in\{0,1\}^s} \cC^{|\omega|} u_{h^\omega}(x) \cdot\cC^{|\omega|+s} g_t(x+h_1^\omega+\cdots+h_s^\omega).
$$
Translating $x$ by $-\sum_\ell h'_\ell$ and using the definition of $U_{h,h'}$, we can rewrite the expression above as
\[
 \E_{h,h'}\E_{x,t}U_{h,h'}(x)\Delta_{h-h'}g_t(x).
\]
Now substitute the definition of $g_t$ and translate $x$ by $P_i(t)$ to rewrite the above as
$$
\E_{h,h'}\E_{x,t} U_{h,h'}(x+P_i(t)) \prod_{j \neq i} \Delta_{h-h'}f_j(x+P_j(t)) = \E_{h,h'} \Lambda_P(\Delta_{h-h'}f_0,\cdots,U_{h,h'},\cdots,\Delta_{h-h'}f_m).
$$
This completes the proof.
\end{proof}

\section{PET induction and concatenation}\label{sec:PET}

In this section we prove Proposition \ref{prop:twisted-pet}. We use Kuca's box norm estimate [ref] for vector polynomials, which combines PET induction and concatenation, with an additional argument ensuring that only lines of the form $L_{ij}$ with $j \in I \cup \{0\}$ appear in the box norm.

\subsection{The unweighted box norm estimate}

We quote the following form of \cite[Proposition~4.3]{Kuca}. Two vector polynomials are called essentially distinct if their difference is nonconstant. Recall that $G = \F_p^k$.

\begin{proposition}\label{prop:pet-box}
Fix positive integers $d,k,m$. Let $Q_1,\ldots,Q_m\in\F_p[t]^k$ be nonconstant, pairwise essentially distinct polynomials of degrees at most $d$. For $0 \leq j < m$, let
\[
 L_j=\langle\lc(Q_m-Q_j) \subset \F_p^k,
\]
where $\langle v\rangle$ denotes the span of a vector $v \in G$, $\lc(Q_m-Q_j)$ denotes the largest coefficient of $Q_m-Q_j$, and we set $Q_0=0$.
Then for all $1$-bounded functions $g_0,\ldots,g_m:G\to\C$,
$$
 \left|\E_{x\in G}\E_{t\in\F_p}
       g_0(x)\prod_{j=1}^m g_j(x+Q_j(t))\right|
 \ll_{d,k,m}\|g_m\|_{L_0^{\times s},\ldots,L_{m-1}^{\times s}}^c+p^{-c},
$$
for some positive integer $s = s(d,k,m)$ and constant $c = c(d,k,m) > 0$.
\end{proposition}

\begin{remark}
It is important in this statement that $Q_1,\cdots,Q_m$ are not required to have zero constant terms, and the implied constants are independent of the coefficients of $Q_1,\cdots, Q_m$.
\end{remark}

In the definition of the counting operator $\Lambda_P$, we shift $t$ by $u^2$ and average over $u \in\F_p$ to write
$$
\Lambda_P(f_0,\cdots,f_m) = \E_{x \in G} \E_{t,u\in\F_p} f_0(x) \prod_{i=1}^m f_i(x+P_i(t+u^2)).
$$
For $t\in \F_p$ and $1 \leq i \leq m$, let $Q_{i,t}(u) = P_i(t+u^2)$, and let $Q_t = (Q_{1,t},\cdots,Q_{m,t})$. Then
\begin{equation}\label{eq:pet-even-average}
\Lambda_P(f_0,\cdots,f_m) = \E_{t \in \F_p} \Lambda_{Q_t}(f_0,\cdots,f_m).
\end{equation}
Note that
\begin{equation}\label{eq:pet-even-leading}
\lc(Q_{i,t}) = \lc(P_i), \ \ \lc(Q_{i,t}-Q_{j,t}) = \lc(P_i-P_j)
\end{equation}
for all $ i \neq j$ and all $t \in \F_p$. In particular, these leading coefficients are independent of $t$.

For each restricted index $i\notin I$, Definition \ref{def:eigenfunctions} gives a frequency $\lambda_{i,z}\in V_i^*$, depending only on $z+V_i$, such that
\[
 f_i(z+v)=f_i(z)e_p(\lambda_{i,z}(v))
\]
for $v \in V_i$. Then
$$
\begin{aligned}
\Lambda_{Q_t}(f_0,\cdots,f_m) &= \E_{x\in G}\E_{u\in\F_p} f_0(x) \prod_{i \in I} f_i(x+Q_{i,t}(u)) \prod_{i\notin I} f_i(x) e_p(\lambda_{i,x}(Q_{i,t}(u))) \\
&= \E_{x \in G} f_0(x)\prod_{i\notin I}f_i(x) \E_{u \in \F_p} \prod_{i \in I} f_i(x+Q_{i,t}(u)) e_p\Big(\sum_{i \notin I} \lambda_{i,x}(Q_{i,t}(u))\Big).
\end{aligned}
$$
The inner sum over $u \in \F_p$ is of the form $\E_u F(u) G(u)$, where
$$
F(u) = F_{t,x}(u) = \prod_{i \in I} f_i(x+Q_{i,t}(u)), \ \ G(u) = e_p\Big(\sum_{i \notin I} \lambda_{i,x}(Q_{i,t}(u))\Big).
$$
Since $G$ is a polynomial phase function, we have $\|G\|_{U^r(\F_p)^*} = 1$ for some $r$ depending on the degrees of $Q_{i,t}$. Thus $|\E_u F(u)G(u)| \leq \|F\|_{U^r(\F_p)}$ and we have
$$
|\Lambda_{Q_t}(f_0,\cdots,f_m)| \leq \E_{x \in G} \|F_{t,x}\|_{U^r(\F_p)}.
$$
By H\"{o}lder's inequality and using the definition of the $U^r$-norm, we have
$$
|\Lambda_{Q_t}(f_0,\cdots,f_m)|^{2^r} \leq \E_{x \in G} \|F_{t,x}\|_{U^r(\F_p)}^{2^r} = \E_{x \in G} \E_{u \in \F_p} \E_{h \in \F_p^r} \prod_{\omega\in\{0,1\}^r} \cC^{|\omega|} F_{t,x}(u+\omega\cdot h).
$$
Plugging in the definition of $F_{t,x}$, we can rewrite the right-hand side above as $\E_{h\in\F_p^r} T_{t,h}$, where
\begin{equation}\label{eq:pet-copied-average}
 T_{t,h}=\E_{x\in G}\E_{u\in\F_p}
   \prod_{i\in I}\prod_{\omega\in\{0,1\}^r}
       \cC^{|\omega|}f_i\bigl(x+Q_{i,t}(u+\omega\cdot h)\bigr).
\end{equation}
Combining the estimates above with \eqref{eq:pet-even-average} and applying H\"{o}lder, we obtain
\begin{equation}\label{eq:pet-weights-removed}
 |\Lambda_P(f_0,\ldots,f_m)|^{2^r}
 \leq\E_t|\Lambda_{Q_t}(f_0,\cdots,f_m)|^{2^r}
 \leq\E_{t,h}|T_{t,h}|.
\end{equation}
Note that only unrestricted indices $i \in I$ remain in the counting operator \eqref{eq:pet-copied-average}.

\subsection{The leading directions and completion of the proof}

Let $\Omega\subset\F_p^r$ consist of those $h$ for which the numbers $\omega\cdot h$ for $\omega\in\{0,1\}^r$ are pairwise distinct. Since $\F_p^r\setminus\Omega$ is the union of at most $2^{2r}$ hyperplanes, we have
$$
 |\F_p^r\setminus\Omega| \leq 2^{2r} p^{r-1}.
$$
In particular, the contribution of $\F_p^r\setminus\Omega$ is bounded by $4^r/p=O(p^{-c})$ since $p$ is sufficiently large. We will focus on $h\in\Omega$ below.

For $h\in\Omega$, we claim that the polynomials
$$
 u\mapsto Q_{i,t}(u+\omega\cdot h)
$$
for $i \in I$ and $\omega \in \{0,1\}^r$
are pairwise essentially distinct. Indeed, consider the difference $u\mapsto Q_{i,t}(u+\omega\cdot h)- Q_{j,t}(u+\omega'\cdot h)$ for some $i,j\in I$ and $\omega,\omega'\in\{0,1\}^r$ with $(i,\omega) \neq (j,\omega')$. If $\omega = \omega'$, then $i \neq j$ and thus the difference is not a constant by \eqref{eq:pet-even-leading}. 
 It remains to consider the case $\omega\neq\omega'$.
If the degrees or the leading coefficients of $Q_{i,t}$ and $Q_{j,t}$ differ, then their difference is again not a constant. In the remaining case, $Q_{i,t}$ and $Q_{j,t}$ have the same even degree $D$ and the same leading coefficient $a \neq 0$, and $\omega\neq\omega'$. Since $Q_{i,t}$ and $Q_{j,t}$ only contain monomials of even degrees, the coefficient of $u^{d-1}$ in the difference is
\begin{equation}\label{eq:pet-even-difference}
 aD(\omega\cdot h-\omega'\cdot h)\ne0,
\end{equation}
since $h \in \Omega$. This proves the claimed essential distinctness.

Thus for each $h \in \Omega$ and $i \in I$, we may apply Proposition \ref{prop:pet-box} to \eqref{eq:pet-copied-average} and obtain
$$
|T_{t,h}| \ll \|f_i\|_{\cB}^c + p^{-c}
$$
for some constant $c>0$, where $\cB$ is a repeated list of lines, each spanned by a vector which is the leading coefficient of $Q_{i,t}$ or of the difference $u\mapsto Q_{i,t}(u) - Q_{j,t}(u+\omega\cdot h)$ for some $(j,\omega) \neq (i,0)$.
We now check that every line appearing in $\cB$ belongs to the collection in the conclusion of Proposition \ref{prop:twisted-pet}.

First, $\lc(Q_{i,t}) = \lc(P_i)$ by \eqref{eq:pet-even-leading}, and thus its span is precisely $L_{i0}$. Now let $(j,\omega) \neq (i,0)$. There are three cases:
\begin{enumerate}
\item If $\omega=0$, then $j\ne i$, and we have $\lc(Q_{i,t}-Q_{j,t}) = \lc(P_i-P_j)$ by \eqref{eq:pet-even-leading}, and hence its span gives the line $L_{ij}$.
\item If $\omega\ne0$ and the degrees or leading coefficients of $Q_{i,t}$ and $Q_{j,t}$ differ, then the leading coefficient of $u\mapsto Q_{i,t}(u) - Q_{j,t}(u+\omega\cdot h)$ is the same as that of $Q_{i,t}-Q_{j,t}$. Hence its span again gives $L_{ij}$.
\item If $\omega \neq 0$ and the two polynomials $Q_{i,t}$ and $Q_{j,t}$ have the same degree $D$ and the same leading coefficient $a=\lc(P_i)$, then by \eqref{eq:pet-even-difference}, the leading coefficient of $u\mapsto Q_{i,t}(u)-Q_{j,t}(u+\omega\cdot h)$ is $-aD(\omega\cdot h)$. Since $h\in\Omega$, this is a nonzero multiple of $a=\lc(P_i)$ and hence it span is $L_{i0}$.
\end{enumerate}

Combining the estimates above with \eqref{eq:pet-weights-removed}, we obtain the desired conclusion
\[
 |\Lambda_P(f_0,\ldots,f_m)|^{2^r}
 \ll\|f_i\|_{\mathcal B}^c+p^{-c}.
\]

\section{Counting with an invariant factor}\label{sec:invariant}

In this section we prove Proposition~\ref{prop:invariant}. The rough strategy is as follows. Suppose that $f_i$ is invariant under translations by a line $L\subset\sum_{i\in I}V_i$. A mixed norm control (Lemma \ref{lem:invariant-energy}, which makes use of the induction hypothesis in Theorem \ref{thm:strengthened}) and an elementary inverse theorem (Lemma \ref{lem:invariant-mixed-inverse}) combine to show that $f_i$ correlates with a structured function $K$ (from the conclusion of Lemma \ref{lem:invariant-mixed-inverse}). In this sketch, let's assume that $f_i$ is such a structured function.
If $L\cap V_i=\{0\}$, then $f_i$ is a product of a $V_i$-invariant function and an $L$-invariant function, so we can increase the number of $L$-invariant functions in the counting operator. If $L\subset V_i$ and $V_i = L\oplus W$, then $f_i$ is a product of an $L$-invariant function $h$ and a $W$-invariant eigenfunction $\chi\in\cE(V_i)$.
We then split $P_i$ into a sum of two components $P_i = Q_i + R_i$, such that $Q_i(t)$ takes values in $L$ and $R_i(t)$ takes values in $W$. Then
$$
f_i(x+P_i(t)) = h(x+P_i(t)) \chi(x+P_i(t)) = h(x+R_i(t)) \chi(x+Q_i(t))
$$
for all $x$ and $t$. This replaces the  polynomial $P_i$ in the original counting operator by $R_i$ which has a smaller coefficient space ($\dim V(R_i) < \dim V(P_i)$), at the cost of adding an polynomial $Q_i$ with a restricted index. This process still decreases the complexity and one can complete the proof by the induction hypothesis of Theorem \ref{thm:strengthened}.

We now turn to the details. Fix an index $i_0 \in I$ such that $f_{i_0}$ is $L$-invariant. Since any $1$-bounded function $f_{i_0}$ can be decomposed as $f_{i_0} = g_{i_0} + h_{i_0}$, where $g_{i_0} = \E(f_{i_0}\vert V_{i_0})$ is $V_{i_0}$-invariant and $h_{i_0}$ satisfies $\E(h_{i_0}\vert V_{i_0}) = 0$, it suffices to prove Proposition \ref{prop:invariant} under the additional assumption that either $f_{i_0}$ is $V_{i_0}$-invariant or $\E(f_{i_0}\vert V_{i_0}) = 0$.

If $f_{i_0}$ is $V_{i_0}$-invariant, then $f_{i_0}(x+P_{i_0}(t)) = f_{i_0}(x)$ and hence
$$
\Lambda_P(f_0,\cdots,f_m) = \E_{x \in G} f_0(x) f_{i_0}(x) \E_{t \in \F_p} \prod_{\substack{1 \leq i \leq m \\i \neq i_0}} f_i(x + P_i(t)).
$$
Since Theorem \ref{thm:strengthened} has already been proved for counting operators with smaller complexity, the right-hand side above is
$$
\E_{x \in G} f_0(x) f_{i_0}(x) \prod_{\substack{1 \leq i \leq m \\ i \neq i_0}} \E(f_i\vert V_i)(x) + O(p^{-c})= M_P(f_0,\cdots,f_m) + O(p^{-c}),
$$
since $f_{i_0} = \E(f_{i_0}\vert V_{i_0})$. Henceforth we assume that $\E(f_{i_0}\vert V_{i_0}) = 0$, and our goal is to prove that
\begin{equation}\label{eq:L-invariant-goal}
\Lambda_P(f_0,\cdots,f_m) \ll p^{-c}.
\end{equation}

\subsection{Two mixed-norm lemmas}

For a subspace $V\subset G$ and a line $L\subset G$, write $\|F\|_{V,L}$ for the box norm associated with the list $(V,L)$. With $\Delta_hF(x)=F(x)\overline{F(x+h)}$, its definition gives
\begin{equation}\label{eq:invariant-mixed-norm}
 \begin{aligned}
 \|F\|_{V,L}^4
 &=\E_{x\in G}\E_{v\in V,h\in L}
 F(x)\overline{F(x+v)}\,\overline{F(x+h)}F(x+v+h)\\
 &=\E_{h\in L}\|\E(\Delta_hF\vert V)\|_2^2 = \E_{h \in L} \|\Delta_hF\|_{U^1(V)}^2.
 \end{aligned}
\end{equation}

\begin{lemma}\label{lem:invariant-energy}
Let the notations and assumptions be as in the statement of Proposition \ref{prop:invariant}. If $f_{i_0}$ is $L$-invariant for some $i_0\in I$, then for every $i\in I\setminus\{i_0\}$ we have
\[
 |\Lambda_P(f_0,\ldots,f_m)|^2
 \leq\|f_i\|_{V_i,L}^2+O(p^{-c}).
\]
\end{lemma}

\begin{proof}
Set $P_0=0$. In the definition of $\Lambda_P(f_0,\cdots,f_m)$, shift $x$ by some $\ell \in L$ and average over all $\ell \in L$, we obtain
\[
 \Lambda_P(f_0,\ldots,f_m)
 =\E_{x,t}f_{i_0}(x+P_{i_0}(t))\E_{\ell\in L}
 \prod_{\substack{0\leq i\leq m\\i\ne i_0}}f_j(x+\ell+P_j(t)),
\]
since $f_{i_0}$ is $L$-invariant.
By Cauchy--Schwarz, we obtain
\[
 |\Lambda_P(f_0,\ldots,f_m)|^2
 \leq\E_{\ell,\ell'\in L}\E_{x,t}
 \prod_{\substack{0\leq i\leq m\\i\ne i_0}}f_j(x+\ell+P_i(t))\overline{f_j(x+\ell'+P_i(t))}.
\]
Replacing $x$ by $x-\ell$ and letting $h = \ell'-\ell$, we obtain
$$
 |\Lambda_P(f_0,\ldots,f_m)|^2 \leq \E_{h \in L} \E_{x,t}  \prod_{\substack{0\leq i\leq m\\i\ne i_0}} \Delta_hf_i(x+P_i(t)).
 $$
For each $h \in L$, the inner average above is the counting operator with $P_{i_0}$ removed from $P$. If $j\notin I$ is a restricted index, then $f_j \in \cE(V_j)$ and hence $\Delta_hf_j \in\cE(V_j)$ as well. Since the complexity $\mathfrak c(P,I)$ strictly decreases when $P_{i_0}$ is removed from $P$ and $i_0$ is removed from $I$, our induction hypothesis implies that 
$$
\E_{x,t}  \prod_{\substack{0\leq i\leq m\\i\ne i_0}} \Delta_hf_i(x+P_i(t)) \leq \min_{i \in I\setminus\{i_0\}} \|\Delta_hf_i\|_{U^1(V_i)} + O(p^{-c}).
$$
Hence, by Cauchy-Scharz and \eqref{eq:invariant-mixed-norm}, for each $i \in I\setminus\{i_0\}$ we have
$$
 |\Lambda_P(f_0,\ldots,f_m)|^4 \leq \E_{h \in L} \|\Delta_hf_i\|_{U^1(V_i)}^2 + O(p^{-c}) = \|f_i\|_{V_i,L}^4 + O(p^{-c}).
 $$
This completes the proof.
\end{proof}

\begin{lemma}\label{lem:invariant-mixed-inverse}
Let $V\subset G$ be a subspace, let $L\subset G$ be a line, and let $F:G\to\C$ be $1$-bounded. There exists a $1$-bounded function $K$ such that
$$
 |\E_x F(x)K(x)|\geq\|F\|_{V,L}^4.
$$
Moreover, $K$ has one of the following forms, with all factors $1$-bounded:
\begin{enumerate}
\item If $L\cap V=\{0\}$, then $K=ah$, where $a$ is $V$-invariant and $h$ is $L$-invariant.
\item If $L\subset V$, choose a complement $W$ with $V=L\oplus W$. Then $K=h\chi$, where $h$ is $L$-invariant and $\chi\in\cE(V)$ is $W$-invariant.
\end{enumerate}
\end{lemma}

\begin{proof}
Suppose first that $L\cap V=\{0\}$. Choose a complement $Z$ with $G=Z\oplus V\oplus L$, and put $F_z(v,b)=F(z+v+b)$ for $z \in Z, v\in V, b \in L$. Then
\[
 \|F\|_{V,L}^4=\E_{z\in Z}\beta_z,\qquad
 \beta_z=\E_{v,v'\in V,b,b'\in L}
 F_z(v,b)\overline{F_z(v',b)}\,\overline{F_z(v,b')}F_z(v',b').
\]
Each $\beta_z$ is nonnegative. By averaging, we may choose $v_z\in V$ and $b_z\in L$ such that
\[
 \mathop{\rm Re}\E_{v,b}
 F_z(v,b)\overline{F_z(v_z,b)}\,\overline{F_z(v,b_z)}F_z(v_z,b_z)
 \geq\beta_z.
\]
Then $K=ah$ where
\[
 a(z+v+b)=\overline{F_z(v_z,b)}F_z(v_z,b_z),
 \qquad h(z+v+b)=\overline{F_z(v,b_z)}
\]
satisfies the desired properties.

Suppose now that $L\subset V$. Choose $Z$ with $G=Z\oplus V$, and write $x=z+w+b$ with $z\in Z$, $w\in W$ and $b\in L$. Define
\[
 \widehat F_{z,w}(\xi)=\E_{b\in L}F(z+w+b)e_p(-\xi(b)),
 \qquad
 \mu_z(\xi)=\E_{w\in W}|\widehat F_{z,w}(\xi)|^2
 \quad(\xi\in L^*).
\]
Fourier orthogonality and Parseval's identity give
\[
 \|F\|_{V,L}^4=\E_{z\in Z}\sum_{\xi\in L^*}\mu_z(\xi)^2,
 \qquad \sum_{\xi\in L^*}\mu_z(\xi)\leq1.
\]
Choose $\xi_z$ maximizing $\mu_z$, and put $\chi(z+w+b)=e_p(-\xi_z(b))$. Thus $\chi\in\cE(V)$ and is $W$-invariant. Choose $h(z+w)$ of modulus at most one so that
\[
 h(z+w)\widehat F_{z,w}(\xi_z)=|\widehat F_{z,w}(\xi_z)|,
\]
and extend $h$ to be $L$-invariant. Since $|\widehat F_{z,w}(\xi_z)|\leq1$, we have
\[
 \E_xF(x)h(x)\chi(x)
 =\E_{z,w}|\widehat F_{z,w}(\xi_z)|
 \geq\E_z\mu_z(\xi_z)
 \geq\E_z\sum_\xi\mu_z(\xi)^2.
\]
This gives the desired function $K = h\chi$.
\end{proof}

\subsection{Proof of Proposition~\ref{prop:invariant}}

We use descending induction on the number of indices $i\in I$
for which $f_i$ is $L$-invariant, starting with the case in which
all these functions are $L$-invariant.

Consider first the base case when $f_i$ is $L$-invariant for all $i \in I$. Choose a projection $T$ on $\sum_{i \in I}V_i$ with kernel $L$, so $v-Tv\in L$ for every $v\in \sum_{i \in I}V_i$. For each $i\in I$, from the $L$-invariance of $f_i$ we have $f_i(x+P_i(t))=f_i(x+TP_i(t))$. Hence
$$
\Lambda_P(f_0,\cdots,f_m) = \E_{x \in G} \E_{t\in\F_p} f_0(x) \prod_{i\in I} f_i(x+TP_i(t)) \prod_{j\notin I} f_j(x+P_j(t)).
$$
Consider the counting operator on the right-hand side above.
Extend $T$ linearly to $G$.
For $i \in I$, the spaces $V(TP_i)=TV_i$ and $W(TP_i)\subset W(P_i)$, since every coordinate of $TP_i$ is a linear combination of the coordinate of $P_i$.
% defined by \eqref{eq:VW}  are precisely $T(V_i)$ and $T(W_i)$. Since $T(W_i) \subset W_i$, 
Thus, the independence condition \eqref{eq:independence} is preserved. 
Moreover,
$$
 \dim\Big(\sum_{i\in I}TV_i\Big)=\dim T\Big(\sum_{i\in I}V_i\Big)=\dim \sum_{i\in I}V_i-1,
$$
because $L\subset \sum_{i\in I}V_i$. Hence this counting operator has smaller complexity. Since Theorem \ref{thm:strengthened} has already been proved for counting operators with smaller complexity, we get
$$
\Lambda_P(f_0,\cdots,f_m) = \E_{x \in G} f_0(x) \prod_{i \in I} \E(f_i\vert T(V_i))(x) \prod_{j \notin I} \E(f_j\vert V_j)(x) + O(p^{-c}).
$$
For $i \in I$, from the $L$-invariance of $f_i$ we have $\E(f_i\vert V_i) = \E(f_i\vert T(V_i))$, and in particular $\E(f_{i_0}\vert T(V_{i_0}))=0$. This proves \eqref{eq:L-invariant-goal} as desired.

We now turn to the induction step. Recall that $f_{i_0}$ is $L$-invariant with $\E(f_{i_0}\vert V_{i_0})=0$. Assume that $f_r$ is not $L$-invariant for some $r \in I$. Define
\[
 d_r(x)=\E_t f_0(x-P_r(t))
 \prod_{\substack{1\leq j\leq m\\j\ne r}}
 f_j(x+P_j(t)-P_r(t)).
\]
Then $d_r$ is $1$-bounded and $\Lambda_P(f_0,\ldots,f_m)=\E_x f_r(x)d_r(x)$. Cauchy--Schwarz gives
\[
 |\Lambda_P(f_0,\ldots,f_m)|^2
 \leq\E_x|d_r(x)|^2
 =\Lambda_P(f_0,\ldots,f_{r-1},\overline{d_r},f_{r+1},\ldots,f_m).
\]
Apply Lemma~\ref{lem:invariant-energy} to the right-hand side above to obtain
$$
 |\Lambda_P(f_0,\ldots,f_m)|^2 \leq \|d_r\|_{V_r,L} + O(p^{-c}).
$$
Then apply Lemma~\ref{lem:invariant-mixed-inverse} to obtain a function $K$ of one of the two forms such that
$$
|\Lambda_P(f_0,\ldots,f_m)|^8 \leq |\E_x d_r(x) K(x)| + O(p^{-c}) = |\Lambda_P(f_0,\cdots,f_{r-1},K,f_{r+1},\cdots,f_m)| + O(p^{-c}).
$$
Hence it suffices to prove that
\begin{equation}\label{eq:L-invariant-K}
|\Lambda_P(f_0,\cdots,f_{r-1},K,f_{r+1},\cdots,f_m)| \ll p^{-c}.
\end{equation}

If $L\cap V_r=\{0\}$, then we have $K=ah$ from the lemma, where $a$ is $V_r$-invariant and $h$ is $L$-invariant. Hence $K(x+P_r(t))=a(x)h(x+P_r(t))$ and we have
$$
\Lambda_P(f_0,\cdots,f_{r-1},K,f_{r+1},\cdots,f_m) = \Lambda_P(f_0a,\cdots,f_{r-1},h,f_{r+1},\cdots,f_m).
$$
Since the counting operator on the right-hand side above has one more $L$-invariant function in $h$, we may apply the induction hypothesis to complete the induction step.

It remains to consider the case $L\subset V_r$. Choose $W$ with $V_r=L\oplus W$, and write $K=h\chi$ as in the second case of Lemma~\ref{lem:invariant-mixed-inverse}, where $h$ is $L$-invariant and $\chi\in\cE(V_r)$ is $W$-invariant. 
Choose projections $S:V_r\rightarrow L$ and $T:V_r\rightarrow W$ such that $v = S(v) + T(v)$ for every $v \in V_r$. Using the decomposition $P_r = SP_r + TP_r$, we can write
$$
K(x+P_r(t)) = h(x+P_r(t)) \chi(x+P_r(t)) = h(x+TP_r(t)) \chi(x+SP_r(t))
$$
by the invariance properties of $h,\chi$. Hence the counting operator in \eqref{eq:L-invariant-K} can be written as
$$
\E_{x \in G} \E_{t \in \F_p} f_0(x) h(x+TP_r(t)) \chi(x + SP_r(t)) \prod_{\substack{1 \leq j \leq m \\j \neq r}} f_j(x+P_j(t)) .
$$

We first note that $V(TP_r)=\text{span}\{Ta_{r,n}:\ n\ge 1\}=T(V_r)=W$, and similarly $V(SP_r)=S(V_r)=L$.
Next one can verify easily that $W(P_r) = W(TP_r)+ W(SP_r)$, which actually is a direct sum using that $\dim W(TP_r)=\dim V(TP_r)=\dim W$ and $\dim W(SP_r)=\dim L$.
%First note that the spaces $W(TP_r)$ and $W(SP_r)$ defined by \eqref{eq:VW} satisfy $W(P_r) = W(TP_r)\oplus W(SP_r)$, and hence the independence assumption \eqref{eq:independence} is satisfied. Next note that the spaces $V(TP_r)$ and $V(SP_r)$ defined by \eqref{eq:VW} are precisely $T(V_r)=W$ and $S(V_r)=L$, respectively. 
Since $\dim W < \dim V_r$ and $\chi \in \cE(V_r) \subset \cE(L)$, this counting operator has smaller complexity, and hence Theorem \ref{thm:strengthened} has already been proved for it by assumption. It is equal to
$$
\E_{x \in G} f_0(x) \E(h\vert W)(x) \E(\chi\vert L)(x) \prod_{\substack{1 \leq j \leq m \\ j\neq r}} \E(f_j\vert V_j)(x) + O(p^{-c}).
$$
Since $i_0 \neq r$ and $\E(f_{i_0}\vert V_{i_0}) = 0$, the expression above is $O(p^{-c})$, proving \eqref{eq:L-invariant-K} as desired.

\section{Obtaining higher-order control along $V_i$}\label{sec:smoothing}

In this section we prove Proposition \ref{prop:high-norm}. 
We recall that we always assume that $m\ge2$.

\begin{lemma}\label{lem:one-transfer}
Let $P = (P_1,\cdots,P_m)$ be vector polynomials in $\F_p[t]^k$ of degree at most $d$ as in the assumption of Theorem \ref{thm:main}, and let $I \subset \{1,\cdots,m\}$. Suppose that the following two statements have already been established for some $i_0 \in I$ and some subspace $L \subset G$:
\begin{enumerate}
\item For any $1$-bounded functions $f_0,\cdots,f_m:G\rightarrow\C$ with $f_j \in \cE(V_j)$ for $j \notin I$, we have
$$
|\Lambda_P(f_0,\cdots,f_m)| \ll_{k,m,d} \|f_{i_0}\|_{\mathcal{B},L}^c + p^{-c}
$$
for some constant $c = c(k,m,d) > 0$,
where $\mathcal{B}$ is a list of $O_{k,m,d}(1)$ subspaces of $G$.
\item For any $1$-bounded functions $f_0,\cdots,f_m:G\rightarrow\C$ with $f_j \in \cE(V_j)$ for $j \notin I$, such that $f_{i_0}$ is $L$-invariant, we have
$$
\Lambda_P(f_0,\cdots,f_m)= M_P(f_0,\cdots,f_m) + O_{k,m,d}(p^{-c})
$$
for some constant $c = c(k,m,d) > 0$.
\end{enumerate}
Then for any $1$-bounded functions $f_0,\cdots,f_m:G\rightarrow\C$ with $f_j \in \cE(V_j)$ for $j \notin I$, we have
$$
 |\Lam_P(f_0,\cdots,f_m)|\ll_{k,m,d} \min_{\substack{1 \leq i \leq m \\i \neq i_0}}\|f_i\|_{\mathcal{B}, V_i}^{c}+p^{-c}
$$
for some constant $c  = c(k,m,d) > 0$.
\end{lemma}

\begin{lemma}\label{lem:one-transfer2}
Let $P = (P_1,\cdots,P_m)$ be vector polynomials in $\F_p[t]^k$ of degree at most $d$ as in the assumption of Theorem \ref{thm:main}, and let $I \subset \{1,\cdots,m\}$. Suppose that the assumptions (1) and (2) in Lemma \ref{lem:one-transfer} have already been established for some $i_0 \notin I$ and $L = V_{i_0}$. Then the conclusion of Lemma \ref{lem:one-transfer} still holds.
\end{lemma}

\subsection{Proof of Lemmas \ref{lem:one-transfer}}

Let $f_0,\cdots,f_m:G\rightarrow\C$ be $1$-bounded functions with $f_j \in \cE(V_j)$ for $j \notin I$, and  let $1 \leq i \leq m$ with $i \neq i_0$. Our goal is to prove that
$$
|\Lambda_P(f_0,\cdots,f_m)| \ll \|f_i\|_{\cB,V_i}^c + p^{-c}.
$$

We first pass to the dual. Define $d_{i_0}: G\rightarrow\C$ by
$$
d_{i_0}(x) = \E_{t \in \F_p} f_0(x-P_{i_0}(t)) \prod_{\substack{1 \leq j \leq m \\ j\neq i_0}} f_j(x+P_j(t)-P_{i_0}(t)).
$$
Then we have
$$
\Lambda_P(f_0,\cdots,f_m) = \E_{x \in G} f_{i_0}(x) d_{i_0}(x) \leq \E_{x \in G} |d_{i_0}(x)|.
$$  
By Cauchy-Schwarz, we have
\begin{equation}\label{eq:one-transfer-dual}
|\Lambda_P(f_0,\cdots,f_m)|^2 \leq \E_{x \in G}|d_i(x)|^2 = \Lambda_P(f_0,\cdots,\overline{d_{i_0}},\cdots,f_m),
\end{equation}
where $\overline{d_{i_0}}$ appears in the $f_{i_0}$-component of the counting operator. By assumption (1) applied to the counting operator $\Lambda_P(f_0,\cdots,\overline{d_{i_0}},\cdots,f_m$ we have
\begin{equation}\label{eq:one-transfer-dual2}
|\Lambda_P(f_0,\cdots,f_m)| \ll \|d_{i_0}\|_{\mathcal{B},L}^c + p^{-c}.
\end{equation}
Write $\mathcal{B} = (H_1,\cdots,H_r)$. Basic properties of the $\|\cdot\|_{\mathcal{B},L}$-norm give
$$
\|d_{i_0}\|_{\mathcal{B},L}^{2^{r+1}} = \E_{h \in H_1\times\cdots\times H_r} \|\E(\Delta_h d_{i_0}\vert L)\|_2^2 = \E_{h \in H_1 \times\cdots\times H_r} \E_{x \in G} \Delta_hd_{i_0}(x) u_h(x),
$$
where $u_h: G\rightarrow\C$ is defined by $u_h(x) = \overline{\E(\Delta_hd_{i_0}\vert L)(x)}$. Let $\delta$ be the expression above. By dual-difference interchange, we have
$$
\E_{h,h' \in H_1\times\cdots\times H_r} |\Lambda_P(\Delta_{h-h'}f_0,\cdots, U_{h,h'},\cdots,\Delta_{h-h'}f_m)| \geq \delta^{2^r}.
$$
Here the function $U_{h,h'}$ appears in the $f_{i_0}$-component of the counting operator $\Lambda_P$ and is defined by
$$
U_{h,h'}(x) = \prod_{\omega \in \{0,1\}^{s-2}} \cC^{|\omega|}u_{h^\omega}\Big(x-\sum_\ell h_{\ell}'\Big),
$$
where $h^\omega_{\ell}=h_\ell$ if $\omega_{\ell}=0$ and $h^\omega_{\ell}=h_\ell'$ if $\omega_{\ell}=1$. Since each $u_{h^\omega}$ is $L$-invariant, $U_{h,h'}$ is also $L$-invariant. 
For $j \notin I$, we have $f_j \in \cE(V_j)$ and hence $\Delta_{h-h'}f_j\in \cE(V_j)$.
 Thus we may apply assumption (2) to obtain
\begin{equation}\label{eq:one-transfer-u}
\Lambda_P(\Delta_{h-h'}f_0,\cdots, U_{h,h'},\cdots,\Delta_{h-h'}f_m) = M_P(\Delta_{h-h'}f_0,\cdots, U_{h,h'},\cdots,\Delta_{h-h'}f_m) + O(p^{-c}).
\end{equation}
Since $i \neq i_0$, the definition of $M_P$ and the triangle inequality implies that
$$
|M_P(\Delta_{h-h'}f_0,\cdots, U_{h,h'},\cdots,\Delta_{h-h'}f_m)| \leq \E_{x \in G} |(\E (\Delta_{h-h'}f_i \vert V_i)(x)|.
$$
It follows that
$$
\delta^{2^r} \leq \E_{h,h' \in H_1\times\cdots\times H_r} \E_{x \in G} |(\E (\Delta_{h-h'}f_i \vert V_i)(x)|+ O(p^{-c}).
$$
The average over $h'$ above can be removed. By Cauchy-Schwarz, we have
$$
\delta^{2^{r+1}} \leq \E_{h\in H_1\times\cdots\times H_r} \E_{x \in G} |(\E (\Delta_{h}f_i \vert V_i)(x)|^2+ O(p^{-c}) = \|f_i\|_{\mathcal{B},V_i}^{2^{r+1}} + O(p^{-c}).
$$
This completes the proof of Lemma \ref{lem:one-transfer}.

\subsection{Proof of Lemma \ref{lem:one-transfer2}}

Let $f_0,\cdots,f_m: G\rightarrow \C$ be $1$-bounded functions with $f_j \in \cE(V_j)$ for $j \notin I$, and let $1 \leq i \leq m$ with $i \neq i_0$. Our goal is to prove that
$$
|\Lambda_P(f_0,\cdots,f_m)| \ll \|f_i\|_{\cB,V_i}^c + p^{-c}.
$$
We proceed as in the proof of Lemma \ref{lem:one-transfer} by defining the dual function $d_{i_0}: G\rightarrow \C$ leading to \eqref{eq:one-transfer-dual}. We cannot apply assumption (1) directly to the counting operator $\Lambda_P(f_0,\cdots,\overline{d_{i_0}},\cdots,f_m)$ since $i_0$ is a restricted index while $d_{i_0}$ may not lie in $\cE(V_{i_0})$.

To remedy this issue, we note from $f_{i_0} \in\cE(V_{i_0})$ that we have
$$
f_{i_0}(x+v) = f_{i_0}(x) e_p(\lambda_x(v))
$$
for all $x \in G$ and $v \in V_{i_0}$, where $\lambda_x \in V_{i_0}^*$ depends only on the coset $x + V_{i_0}$ in $G/V_{i_0}$ (see \eqref{eq:eigenfn-property}). Since
$$
\Lambda_P(f_0,\cdots,f_m) = \E_{x \in G} f_{i_0}(x) d_{i_0}(x),
$$
by shifting the variable $x$ by $v \in V_{i_0}$ and averaging over $v$ we obtain
\begin{equation}\label{eq:one-transfer2-g}
\Lambda_P(f_0,\cdots,f_m) = \E_{x \in G} \E_{v \in V_{i_0}} f_{i_0}(x) e_p(\lambda_x(v)) d_{i_0}(x+v) = \E_{x \in G} f_{i_0}(x) g(x),
\end{equation}
where $g: G\rightarrow \C$ is defined by
$$
g(x) = \E_{v \in V_{i_0}} d_{i_0}(x+v) e_p(\lambda_x(v)).
$$
Applying \eqref{eq:one-transfer2-g} with all $f_{i_0} = \overline{g}$, we also obtain
$$
\Lambda_P(f_0,\cdots,\overline{g},\cdots,f_m) = \E_{x \in G} |g(x)|^2,
$$
where $\overline{g}$ appears in the $f_{i_0}$-component of the counting operator.
Hence it follows by Cauchy-Schwarz that
$$
|\Lambda_P(f_0,\cdots,f_m)|^2 \leq \E_{x \in G}|g(x)|^2 = \Lambda_P(f_0,\cdots,\overline{g},\cdots,f_m).
$$
By definition, we have $g \in \cE(V_{i_0})$, and hence by assumption (1) we have
$$
\Lambda_P(f_0,\cdots,\overline{g},\cdots,f_m) \ll \|g\|_{\cB,L}^c + p^{-c}.
$$
Combining this with $\|g\|_{\cB,L} \leq \|d_{i_0}\|_{\cB,L}$ from Lemma~\ref{lem:frequency-projection}, we arrive at \eqref{eq:one-transfer-dual2}. Now proceed exactly as in the proof of Lemma \ref{lem:one-transfer}. Note that we may still apply assumption (2) to obtain \eqref{eq:one-transfer-u} since the $L$-invariance of $U_{h,h'}$ and the assumption that $L = V_{i_0}$ implies that $U_{h,h'} \in \cE(V_{i_0})$.

\subsection{Deduction of Proposition \ref{prop:high-norm}}

Let $k \in I$ be arbitrary. By Proposition \ref{prop:twisted-pet}, we have
\begin{equation}\label{eq:norm-smoothing}
|\Lambda_P(f_0,\cdots,f_m)| \ll \|f_k\|_{\cB_k}^c + p^{-c},
\end{equation}
where $\cB_k$ is a (repeated) list of $O(1)$ lines from $\{L_{kj}\}$ with $j \in I \cup \{0\}$ as in the statement of Proposition \ref{prop:twisted-pet}. Note that all these lines $L_{kj}$ lie in $\sum_{i \in I}V_i$ by definition. We will replace all lines in $\cB_k$ by $V_k$ one by one, and each replacement will involve two applications of Lemma \ref{lem:one-transfer} or \ref{lem:one-transfer2}.

First write $\cB_k = (\cB, L_k)$ for some line $L_k \subset \sum_{i\in I}V_i$. We will apply Lemma \ref{lem:one-transfer} with $i_0 = k$ and $L = L_k$. Assumption (1) is satisfied by \eqref{eq:norm-smoothing} and assumption (2) is satisfied by Proposition \ref{prop:invariant} (making use of the hypothesis that Theorem \ref{thm:strengthened} has already been established for counting operators with smaller complexity). Thus Lemma \ref{lem:one-transfer}  implies that
\begin{equation}\label{eq:norm-smoothing2}
|\Lambda_P(f_0,\cdots,f_m)| \ll \|f_\ell\|_{\mathcal{B}, V_\ell}^c + p^{-c}
\end{equation}
for any $1$-bounded functions $f_0,\cdots,f_m:G\rightarrow\C$ with $f_j \in \cE(V_j)$ for $j \notin I$, and any $1 \leq \ell \leq m$ with $\ell \neq k$. Pick an arbitrary such index $\ell$.

We will now apply Lemma \ref{lem:one-transfer} (when $\ell \in I$) or Lemma \ref{lem:one-transfer2} (when $\ell\notin I$) with $i_0=\ell$ and $L = V_\ell$. Assumption (1) is satisfied by \eqref{eq:norm-smoothing2}. To establish assumption (2) when $f_\ell$ is $V_\ell$-invariant, note that in this case we have
$$
\Lambda_P(f_0,\cdots,f_m) = \E_{x \in G} (f_0f_\ell)(x) \E_{t \in \F_p} \prod_{\substack{1 \leq j \leq m \\j \neq \ell}} f_j(x+P_j(t)).
$$
Thus assumption (2) follows since Theorem \ref{thm:strengthened} has already been established for counting operators with lower complexity. Thus Lemma \ref{lem:one-transfer} or Lemma \ref{lem:one-transfer2} implies that
$$ 
|\Lambda_P(f_0,\cdots,f_m)| \ll \|f_k\|_{\mathcal{B}, V_k}^c + p^{-c}
$$
for any $1$-bounded functions $f_0,\cdots,f_m:G\rightarrow\C$ with $f_j \in \cE(V_j)$ for $j \notin I$. This replaces one of the lines in $\cB_k$ from \eqref{eq:norm-smoothing} by $V_k$. Repeating this process to each of the other lines in $\cB_k$ completes the proof.

\section{Degree lowering and completion of proof}\label{sec:degree-lowering}

In this section we prove Theorem \ref{thm:strengthened}, which implies Theorem \ref{thm:main}, and deduce Corollary \ref{cor:szemeredi}.

\subsection{Proof setup}

Let $P = (P_1,\cdots,P_m)$ be as in the assumption, and let $I \subset \{1,\cdots,m\}$  be a subset. Let $f_0,\cdots,f_m: G\rightarrow\C$ be $1$-bounded functions with $f_j \in \cE(V_j)$ for $j \notin I$.  We induct on the complexity $\mathfrak c(P,I)$.

If $I = \emptyset$ then we are done by Lemma~\ref{lem:all-restricted}. Henceforth assume that $I$ is nonempty and fix an index $i \in I$. It suffices to show that if  $\E(f_k\vert V_k) = 0$ for some $1 \leq k\leq m$ with $k \neq i$ then
\begin{equation}\label{eq:degree-lowering-goal}
|\Lambda_P(f_0,\cdots,f_m)| \ll p^{-c}.
\end{equation}
Indeed, once \eqref{eq:degree-lowering-goal} is proven, Theorem \ref{thm:strengthened} can be deduced as follows. For each $1 \leq j \leq m$ with $j \neq i$, make the decomposition $f_j = g_j + h_j$ where $g_j = \E(f_j\vert V_j)$, so that $\E(h_j\vert V_j) = 0$. Then $\Lambda_P(f_0,\cdots,f_m)$ can be decomposed into $2^{m-1}$ terms. Applying \eqref{eq:degree-lowering-goal} to each of the $2^{m-1}-1$ terms involving at least one $h_j$, we obtain
$$
\Lambda_P(f_0,\cdots,f_m) = \E_{x \in G} f_0(x) \E_{t \in \F_p}f_i(x+P_i(t)) \prod_{\substack{1 \leq j \leq m \\ j\neq i}} g_j(x+P_j(t)) + O(p^{-c}).
$$
Since $g_j$ is constant on each coset of $V_j$, we have $g_j(x + P_j(t)) = g_j(x)$ and hence
$$
\Lambda_P(f_0,\cdots,f_m) = \E_{x \in G} f_0(x)\prod_{\substack{1 \leq j \leq m \\ j \neq i}}g_j(x) \E_{t \in \F_p} f_i(x+P_i(t)) + O(p^{-c}).
$$
By Lemma~\ref{lem:one-shift}, the inner average over $t$ can be replaced by $\E(f_i\vert V_i)(x)$ at the cost of an error $E(x)$ satisfying $\E_{x \in G} |E(x)|^2 \ll p^{-1}$. Since $g_j = \E(f_j\vert V_j)$, this leads to the conclusion of Theorem \ref{thm:strengthened}.

It remains to establish \eqref{eq:degree-lowering-goal}, under the assumption that $\E(f_k\vert V_k)=0$ for some $1 \leq k\leq m$ with $k \neq i$. Fix such an index $k$. Suppose, for the purpose of contradiction, that \eqref{eq:degree-lowering-goal} fails. Consider the dual function $d_i:G\rightarrow \C$ by
$$
d_i(x) = \E_{t \in \F_p} f_0(x-P_i(t)) \prod_{\substack{1 \leq j \leq m \\j \neq i}} f_j(x+P_j(t)-P_i(t)).
$$
Then we have
$$
p^{-c} \leq |\Lambda_P(f_0,\cdots,f_m)| = |\E_{x \in G} f_i(x) d_i(x)| \leq \E_{x \in G} |d_i(x)|.
$$
By Cauchy-Schwarz, we have
$$
p^{-2c} \leq \E_{x \in G} |d_i(x)|^2 = \E_{x \in G} \overline{d_i(x)}  \E_{t \in \F_p} f_0(x-P_i(t)) \prod_{\substack{1 \leq j \leq m \\j \neq i}} f_j(x+P_j(t)-P_i(t)).
$$
After shifting $x$ by $P_i(t)$, the right-hand side above is precisely $\Lambda_P(f_0,\cdots,\overline{d_i},\cdots,f_m)$, where $\overline{d_i}$ replaces $f_i$ compared to the original counting operator. Thus by the higher-order control (Proposition \ref{prop:high-norm}) making use of the induction hypothesis, we can deduce that $\|d_i\|_{U^s(V_i)} \gg p^{-c}$.

\begin{lemma}[Degree lowering]\label{lem:degree-lowering}
For any $s \geq 3$ we have
$$
\|d_i\|_{U^s(V_i)}^{2^{2s-2}} \ll \|d_i\|_{U^{s-1}(V_i)} + p^{-c}.
$$
\end{lemma}

In the rest of this section we write $F = d_i$ and $V = V_i$. Once Lemma \ref{lem:degree-lowering} is proved, we have $\|d_i\|_{U^2(V_i)} \gg p^{-c}$ and we can deduce \eqref{eq:degree-lowering-goal} as follows. Since
$$
\|F\|_{U^2(V)}^4 = \E_{z \in Z} \sum_{\xi\in V^*} |\wh{F_z}(\xi)|^4 \leq \E_{z \in Z} \max_{\xi \in V^*}  |\wh{F_z}(\xi)| = \sup_{\chi \in \cE(V)} \Big|\E_{x \in G} F(x)\chi(x)\Big|,
$$
it follows that there exists $\chi \in \cE(V)$ such that
$$
p^{-c} \leq \E_{x \in G} F(x)\chi(x) = \Lambda_P(f_0, \cdots, \chi, \cdots, f_m),
$$
where $\chi$ appears in the $f_i$-component of the counting operator. Note that the counting operator $\Lambda_P(f_0,\cdots,\chi,\cdots,f_m)$ has lower complexity since it moves the unrestricted index $i \in I$ to restricted, \eqref{eq:degree-lowering-goal} follows from the induction hypothesis by using the assumption $\E(f_k\vert V_k)=0$.

\subsection{Proof of Lemma \ref{lem:degree-lowering}}

Let $\delta = \|F\|_{U^s(V)}$. We may assume that $\delta > p^{-c}$.
Choose representatives for the cosets of $V$ to form $Z \subset G$ such that $G = V\oplus Z$. For $z \in Z$, write $F_z: V\rightarrow\C$ for the function defined by $F_z(v) = F(z+v)$. By basic properties of the $U^s(V)$-norm, we have
$$
\E_{(z,h) \in Z\times V^{s-2}} \sum_{\xi \in V^*} |\wh{\Delta_hF_z}(\xi)|^4 = \delta^{2^s}.
$$
For $(z,h) \in Z \times V^{s-2}$, choose $\phi_{z,h} \in V^*$ such that
$$
|\wh{\Delta_hF_z}(\phi_{z,h})| = \max_{\xi \in V^*}|\wh{\Delta_hF_z}(\xi)|,
$$
and then choose $\zeta_{z,h} \in \C$ with $|\zeta_{z,h}|=1$ such that
$$
|\wh{\Delta_hF_z}(\phi_{z,h})| = \zeta_{z,h} \wh{\Delta_hF_z}(\phi_{z,h}).
$$
Then we have
\begin{equation}\label{eq:c(z,h)}
c(z,h) := \max_{\xi \in V^*}|\wh{\Delta_hF_z}(\xi)| = \zeta_{z,h}\E_{v \in V} \Delta_hF_z(v)e_p(\phi_{z,h}(v))
\end{equation}
and
$$
\delta^{2^s} \leq \E_{(z,h) \in Z \times V^{s-2}} c(z,h)^2 \sum_{\xi \in V^*} |\wh{\Delta_hF_z}(\xi)|^2 \leq \E_{(z,h) \in Z \times V^{s-2}}c(z,h)^2.
$$
For $h \in V^{s-2}$, define $u_h: G\rightarrow\C$ by 
$$
u_h(z+v) = c(z,h)\zeta_{z,h}e_p(\phi_{z,h}(v))
$$
for $z \in Z$ and $v \in V$. Then $u_h \in \cE(V)$ and 
$$
\begin{aligned}
\E_{x \in G} \E_{h \in V^{s-2}} \Delta_hF(x) u_h(x) &= \E_{(z,h) \in Z\times V^{s-2}} \E_{v \in V} \Delta_hF(z+v) u_h(z+v) \\
&= \E_{(z,h) \in Z\times V^{s-2}} c(z,h)\zeta_{z,h} \E_{v \in V} \Delta_hF_z(v)e_p(\phi_{z,h}(v)) \\
&= \E_{(z,h) \in Z\times V^{s-2}} c(z,h)^2 \geq \delta^{2^s}.
\end{aligned}
$$
By the dual-difference interchange (Lemma \ref{lem:dual-difference}), we have
\begin{equation}\label{eq:degree-lowering-ddchange}
\E_{h,h' \in V^{s-2}} |\Lambda_P(\Delta_{h-h'}f_0,\cdots, U_{h,h'}, \cdots,\Delta_{h-h'}f_m)| \geq \delta^{2^s\cdot 2^{s-2}}.
\end{equation}
Here the function $U_{h,h'}$ appears in the $f_i$-component of the counting operator $\Lambda_P$ and is defined by
$$
U_{h,h'}(x) = \prod_{\omega \in \{0,1\}^{s-2}} \cC^{|\omega|}u_{h^\omega}\Big(x-\sum_\ell h_{\ell}'\Big),
$$
where $h^\omega_{\ell}=h_\ell$ if $\omega_{\ell}=0$ and $h^\omega_{\ell}=h_\ell'$ if $\omega_{\ell}=1$. Since each $u_{h^\omega} \in \cE(V)$, we have $U_{h,h'} \in \cE(V)$. 
For $j \notin I$, we have $f_j \in \cE(V)$ and hence $\Delta_{h-h'}f_j\in \cE(V)$. 
Hence the counting operator in \eqref{eq:degree-lowering-ddchange} is of lower complexity since the index $i$ becomes a restricted index. It follows by induction hypothesis that
$$
\E_{h,h' \in V^{s-2}} |M_P(\Delta_{h-h'}f_0,\cdots, U_{h,h'}, \cdots,\Delta_{h-h'}f_m)| \gg \delta^{2^s\cdot 2^{s-2}},
$$
since $\delta > p^{-c}$ and $c$ is sufficiently small.
By the definition of $M_P$ and isolating the term from $U_{h,h'}$, we obtain
$$
\E_{h,h' \in V^{s-2}} \E_{x \in G} |\E(U_{h,h'}\vert V)(x)| \gg \delta^{2^s\cdot 2^{s-2}}.
$$
From the formula for $U_{h,h'}$ and the definition of $u_{h^\omega}$ we have
$$
U_{h,h'}(z+v) = \prod_{\omega \in \{0,1\}^{s-2}} \cC^{|\omega|} c(z,h^\omega) \zeta_{z,h^\omega} e_p\Big(\phi_{z,h^\omega}\Big(v-\sum_\ell h_\ell'\Big)\Big).
$$
for $z \in Z$ and $v \in V$. For fixed $h,h',z$, the expression above is a linear phase of $v$ with frequency
$$
\Phi_{z,h,h'} = \sum_{\omega \in \{0,1\}^{s-2}} (-1)^\omega \phi_{z,h^\omega}.
$$
It follows that
$$
|\E(U_{h,h'}\vert V)(z)| = \Big|\E_{v \in V} U_{h,h'}(z+v)\Big| = \begin{cases} \prod_{\omega \in \{0,1\}^{s-2}} c(z,h^\omega), & \text{if }\Phi_{z,h,h'}=0, \\
0, & \text{if }\Phi_{z,h,h'} \neq 0.\end{cases}
$$
Therefore we have
$$
\delta^{2^s\cdot 2^{s-2}} \ll \E_{z\in Z}\E_{h,h'\in V^{s-2}} 1_{\{\Phi_{z,h,h'}=0\}} \prod_{\omega \in \{0,1\}^{s-2}} c(z,h^\omega).
$$
For $h' \in V^{s-2}$, define the weight
$$
w_{h'}(z,h) = \zeta_{z,h}  1_{\{\Phi_{z,h,h'}=0\}} \prod_{\substack{\omega \in \{0,1\}^{s-2} \\ \omega\neq 0}} c(z,h^\omega)
$$
for $(z,h) \in Z \times V^{s-2}$. By \eqref{eq:c(z,h)}, for $h' \in V^{s-2}$ we have
$$
\E_{(z,h) \in Z \times V^{s-2}} w_{h'}(z,h) \E_{v \in V}\Delta_hF_z(v)e_p(\phi_{z,h}(v)) = \E_{(z,h) \in Z \times V^{s-2}}  1_{\{\Phi_{z,h,h'}=0\}} \prod_{\omega \in \{0,1\}^{s-2}} c(z,h^\omega).
$$
Thus there exists $h' \in V^{s-2}$ such that
$$
\E_{(z,h) \in Z \times V^{s-2}} w_{h'}(z,h) \E_{v \in V}\Delta_hF_z(v)e_p(\phi_{z,h}(v)) \gg \delta^{2^s\cdot 2^{s-2}}.
$$
Fix such $h'$ from now on, and denote by $J$ the left-hand side above. Decompose the sum in $\Phi_{z,h,h'}$ according to the position of the first nonzero coordinate of $\omega$ as follows. For $1 \leq \ell \leq s-2$, let
$$
\psi_{\ell,z,h,h'} = \sum_{\substack{\omega \in \{0,1\}^{s-2} \\ \omega_{\ell}=1,\omega_1=\cdots=\omega_{\ell-1}=0}} (-1)^{|\omega|+1}\phi_{z,h^\omega}.
$$
Then the condition $\Phi_{z,h,h'}=0$ is equivalent to
$$
\phi_{z,h} = \sum_{1\leq\ell\leq s-2} \psi_{\ell,z,h,h'}.
$$
Since $w_{h'}(z,h) =0$ unless $\Phi_{z,h,h'}=0$, we have
$$
J = \E_{(z,h) \in Z \times V^{s-2}} w_{h'}(z,h) \E_{v \in V} \Delta_h F_z(v) e_p\Big(\sum_\ell \psi_{\ell,z,h,h'}(v)\Big).
$$
By Cauchy-Schwarz, we have
$$
J^2 \leq \E_{(z,h) \in Z \times V^{s-2}} \Big| \E_{v \in V} \Delta_h F_z(v) e_p\Big(\sum_\ell \psi_{\ell,z,h,h'}(v)\Big)\Big|^2.
$$
Expanding the square and making a change of variables, we have
$$
\begin{aligned}
J^2 &\leq \E_{(z,h) \in Z \times V^{s-2}} \E_{v,u \in V} \Delta_hF_z(v+u)\overline{\Delta_hF_z(v)} e_p\Big(\sum_\ell \psi_{\ell,z,h,h'}(v+u)-\psi_{\ell,z,h,h'}(v)\Big) \\
&= \E_{(z,h) \in Z \times V^{s-2}} \E_{v,u \in V}  \Delta_h\Delta_uF_z(v) e\Big(\sum_{\ell} \psi_{\ell,z,h,h'}(u)\Big) \\
&\leq \E_{z \in Z, u \in V} \Big|\E_{v\in V,h \in V^{s-2}} \Delta_h(\Delta_uF_z)(v) \prod_{\ell} e(\psi_{\ell,z,h,h'}(u))\Big|.
\end{aligned}
$$
For fixed $z,u$, note that the function $(v,h)\mapsto \psi_{\ell,z,h,h'}(u)$ does not dependent on the $\ell$th component of $h$. It follows that the inner average above is bounded by $\|\Delta_uF_z\|_{U^{s-2}(V)}$, and hence
$$
\begin{aligned}
\delta^{2^{2s}} &\ll J^2 \leq \E_{z \in Z, u\in V} \|\Delta_uF_z\|_{U^{s-2}(V)} \leq \Big(\E_{z \in Z, u\in V} \|\Delta_uF_z\|_{U^{s-2}(V)}^{2^{s-2}}\Big)^{1/2^{s-2}} \\
&= \Big(\E_{z \in Z} \|F_z\|_{U^{s-1}(V)}^{2^{s-1}}\Big)^{1/2^{s-2}} = \|F\|_{U^{s-1}(V)}^2.
\end{aligned}
$$
This completes the proof.


\begin{thebibliography}{99}
\bibitem{Bergelson} V. Bergelson, \emph{Weakly mixing PET}, Ergodic Theory Dynam. Systems \textbf{7} (1987), no.~3, 337--349. \href{https://doi.org/10.1017/S0143385700004090}{doi:10.1017/S0143385700004090}.


\bibitem{BHMP}
V. Bergelson, B. Host, R. McCutcheon, and F. Parreau,
\emph{Aspects of uniformity in recurrence},
Colloq. Math. \textbf{84/85} (2000), part~2, 549--576.
\href{https://doi.org/10.4064/cm-84/85-2-549-576}{doi:10.4064/cm-84/85-2-549-576}.



\bibitem{BL} V. Bergelson and A. Leibman, \emph{Polynomial extensions of van der Waerden's and Szemer\'edi's theorems}, J. Amer. Math. Soc. \textbf{9} (1996), no.~3, 725--753. \href{https://doi.org/10.1090/S0894-0347-96-00194-4}{doi:10.1090/S0894-0347-96-00194-4}.
\bibitem{BLM} V. Bergelson, A. Leibman and R. McCutcheon, \emph{Polynomial Szemer\'edi theorems for countable modules over integral domains and finite fields}, J. Anal. Math. \textbf{95} (2005), 243--296. \href{https://doi.org/10.1007/BF02791504}{doi:10.1007/BF02791504}.

\bibitem{BC} J. Bourgain and M.-C. Chang, \emph{Nonlinear Roth type theorems in finite fields}, Israel J. Math. \textbf{221} (2017), no.~2, 853--867. \href{https://doi.org/10.1007/s11856-017-1577-9}{doi:10.1007/s11856-017-1577-9}.

\bibitem{DLS} D. Dong, X. Li and W. Sawin, \emph{Improved estimates for polynomial Roth type theorems in finite fields}, J. Anal. Math. \textbf{141} (2020), no.~2, 689--705. \href{https://doi.org/10.1007/s11854-020-0113-8}{doi:10.1007/s11854-020-0113-8}.

\bibitem{DFKMS} S. Donoso, A. Ferr\'e Moragues, A. Koutsogiannis and W. Sun, \emph{Decomposition of multicorrelation sequences and joint ergodicity}, Ergodic Theory Dynam. Systems \textbf{44} (2024), no.~2, 432--480. \href{https://doi.org/10.1017/etds.2023.30}{doi:10.1017/etds.2023.30}.

\bibitem{FK} N. Frantzikinakis and B. Kuca, \emph{Joint ergodicity for commuting transformations and applications to polynomial sequences}, Invent. Math. \textbf{239} (2025), 621--706. \href{https://doi.org/10.1007/s00222-024-01313-w}{doi:10.1007/s00222-024-01313-w}. \href{https://arxiv.org/abs/2207.12288v2}{arXiv:2207.12288v2}.

\bibitem{GMZ}
J. Guo, C. Miao, and G. Zhan,
\emph{A multidimensional Szemer\'edi theorem in integers},
preprint (2026).
\href{https://arxiv.org/abs/2605.06360v1}{arXiv:2605.06360v1}.



\bibitem{HLY} R. Han, M.~T. Lacey and F. Yang, \emph{A polynomial Roth theorem for corners in finite fields}, Mathematika \textbf{67} (2021), no.~4, 885--896. \href{https://doi.org/10.1112/mtk.12108}{doi:10.1112/mtk.12108}.

\bibitem{HZ}
D. He and X. Zhu,
\emph{A Roth theorem in $\mathbb{R}^2$ and a related ergodic theorem},
preprint (2026).
\href{https://arxiv.org/abs/2607.05124v1}{arXiv:2607.05124v1}.



\bibitem{Hong} G.-D. Hong, \emph{On weighted multilinear polynomial averages in finite fields}, preprint (2025). \href{https://arxiv.org/abs/2507.14414v1}{arXiv:2507.14414v1}.
\bibitem{HongSim} G.-D. Hong, \emph{Simultaneous popular polynomial differences in higher dimensions}, preprint (2026). \href{https://arxiv.org/abs/2610.30335v1}{arXiv:2610.30335v1}.

\bibitem{KMPWW}
D. Kosz, M. Mirek, S. Peluse, R. Wan, and J. Wright,
\emph{The multilinear circle method and a question of Bergelson},
Ann. of Math. (2), to appear.
\href{https://arxiv.org/abs/2411.09478v4}{arXiv:2411.09478v4}.


\bibitem{KKS} A. Koutsogiannis, B. Kuca and W. Sun, \emph{Structure of 2-step nilpotent ergodic averages for distinct-degree polynomials}, preprint (2026). \href{https://arxiv.org/abs/2607.29368v2}{arXiv:2607.29368v2}.
\bibitem{KKL} N. Kravitz, B. Kuca and J. Leng, \emph{Quantitative concatenation for polynomial box norms}, Adv. Math. \textbf{489} (2026), 110820. \href{https://doi.org/10.1016/j.aim.2026.110820}{doi:10.1016/j.aim.2026.110820}.
\bibitem{KucaDistinct} B. Kuca, \emph{Multidimensional polynomial Szemer\'edi theorem in finite fields for polynomials of distinct degrees}, Israel J. Math. \textbf{259} (2024), no.~2, 589--620. \href{https://doi.org/10.1007/s11856-023-2551-3}{doi:10.1007/s11856-023-2551-3}.
\bibitem{Kuca} B. Kuca, \emph{Multidimensional polynomial patterns over finite fields: Bounds, counting estimates and Gowers norm control}, Adv. Math. \textbf{448} (2024), 109700. \href{https://doi.org/10.1016/j.aim.2024.109700}{doi:10.1016/j.aim.2024.109700}. \href{https://arxiv.org/html/2304.10793v2}{arXiv:2304.10793v2}.
\bibitem{LSTW} J. Li, C.-Y. Shen, T.~T. Truong and L. Wu, \emph{Polynomial corners in finite fields beyond the distinct-degree case}, preprint (2026). \href{https://arxiv.org/abs/2609.13891v1}{arXiv:2609.13891v1}.
\bibitem{Lim} Z.~L. Lim, \emph{Uniform nonlinear Szemer\'edi theorem for corners in finite fields}, preprint (2025). \href{https://arxiv.org/abs/2501.04887v1}{arXiv:2501.04887v1}.
\bibitem{PeluseThree} S. Peluse, \emph{Three-term polynomial progressions in subsets of finite fields}, Israel J. Math. \textbf{228} (2018), no.~1, 379--405. \href{https://doi.org/10.1007/s11856-018-1768-z}{doi:10.1007/s11856-018-1768-z}.
\bibitem{Peluse} S. Peluse, \emph{On the polynomial Szemer\'edi theorem in finite fields}, Duke Math. J. \textbf{168} (2019), no.~5, 749--774. \href{https://doi.org/10.1215/00127094-2018-0051}{doi:10.1215/00127094-2018-0051}.
\bibitem{PeluseBounds} S. Peluse, \emph{Bounds for sets with no polynomial progressions}, Forum Math. Pi \textbf{8} (2020), e16, 55~pp. \href{https://doi.org/10.1017/fmp.2020.11}{doi:10.1017/fmp.2020.11}.
\bibitem{PP} S. Peluse and S. Prendiville, \emph{Quantitative bounds in the nonlinear Roth theorem}, Invent. Math. \textbf{238} (2024), no.~3, 865--903. \href{https://doi.org/10.1007/s00222-024-01293-x}{doi:10.1007/s00222-024-01293-x}.

\bibitem{PPS}
S. Peluse, S. Prendiville, and X. Shao,
\emph{Bounds in a popular multidimensional nonlinear Roth theorem},
J. London Math. Soc. (2) \textbf{110} (2024), no.~5, e70019.
\href{https://doi.org/10.1112/jlms.70019}{doi:10.1112/jlms.70019}.



\bibitem{TZ} T. Tao and T. Ziegler, \emph{Concatenation theorems for anti-Gowers-uniform functions and Host--Kra characteristic factors}, Discrete Anal. (2016), Paper No.~13. \href{https://doi.org/10.19086/da.850}{doi:10.19086/da.850}.
\bibitem{Weil} A. Weil, \href{https://doi.org/10.1073/pnas.34.5.204}{\emph{On some exponential sums}}, Proc. Natl. Acad. Sci. U.S.A. \textbf{34} (1948), no.~5, 204--207.
\end{thebibliography}
\end{document}